\documentclass[10pt,a4paper]{article}

\usepackage[T1]{fontenc}
\usepackage{lmodern}
\usepackage[a4paper,margin=25mm]{geometry}
\usepackage{amsmath,amssymb,amsthm,mathtools,bm,mathrsfs}
\usepackage{booktabs,array,tabularx,enumitem,float}
\usepackage{microtype,needspace}
\usepackage{titlesec}
\usepackage[numbers,sort&compress]{natbib}
\usepackage[hidelinks]{hyperref}
\hypersetup{
  pdftitle={Morse index, nullity and Jacobi fields of homogeneous minimal two-spheres in complex hyperquadrics},
  pdfauthor={Qing Chen, Shuping Huang and Shoulan Yi},
  pdfsubject={Morse index, nullity, complex hyperquadrics, homogeneous minimal two-spheres},
  pdfkeywords={Morse index, normal nullity, homogeneous minimal two-sphere, complex hyperquadric, Jacobi fields}
}
\allowdisplaybreaks[2]
\setlist{topsep=4pt,itemsep=2.5pt,parsep=1pt,leftmargin=*}
\titleformat{\section}{\normalfont\large\bfseries}{\thesection}{1em}{}
\titleformat{\subsection}{\normalfont\normalsize\bfseries}{\thesubsection}{1em}{}
\titleformat{\paragraph}[runin]{\normalfont\normalsize\bfseries}{}{0pt}{}[.]
\titlespacing{\section}{0pt}{15pt plus 3pt minus 2pt}{7pt plus 2pt minus 1pt}
\titlespacing{\subsection}{0pt}{11pt plus 2pt minus 1pt}{5pt plus 1pt minus 1pt}
\titlespacing*{\paragraph}{0pt}{5pt plus 1pt minus 1pt}{0.7em}
\numberwithin{equation}{section}
\theoremstyle{plain}
\newtheorem{theorem}{Theorem}[section]
\newtheorem{proposition}[theorem]{Proposition}
\newtheorem{lemma}[theorem]{Lemma}

\theoremstyle{definition}

\theoremstyle{remark}
\newtheorem{remark}[theorem]{Remark}
\newtheorem*{remark*}{Remark}
\newcommand{\R}{\mathbb R}
\newcommand{\C}{\mathbb C}
\newcommand{\CP}{\mathbb CP}
\newcommand{\SU}{\mathrm{SU}}

\newcommand{\Ind}{\operatorname{Ind}}
\newcommand{\Nul}{\operatorname{Nul}}
\newcommand{\Hom}{\operatorname{Hom}}

\newcommand{\dd}{\mathrm d}

\begin{document}

\begin{center}
{\fontsize{17.28}{21}\selectfont\bfseries
Morse index, nullity and Jacobi fields of homogeneous minimal two-spheres in complex hyperquadrics\par}
\vspace{14pt}
{\large Qing Chen\textsuperscript{1,2}, Shuping Huang\textsuperscript{3,*}
\& Shoulan Yi\textsuperscript{3}\par}
\vspace{12pt}
{\footnotesize\itshape
\textsuperscript{1}School of Mathematical Sciences and Institute of Geometry and Physics,
University of Science and Technology of China, 96 Jinzhai Road,
Hefei 230026, Anhui Province, China;\par
\textsuperscript{2}Hefei National Laboratory, Hefei 230088, China;\par
\textsuperscript{3}School of Mathematical Sciences,
University of Science and Technology of China,
Hefei, Anhui 230026, P. R. China\par
\vspace{4pt}
Email: \href{mailto:qchen@ustc.edu.cn}{qchen@ustc.edu.cn},
\href{mailto:hsp@mail.ustc.edu.cn}{hsp@mail.ustc.edu.cn},
\href{mailto:yishoulan@mail.ustc.edu.cn}{yishoulan@mail.ustc.edu.cn}\par}
\end{center}
\begingroup
\renewcommand{\thefootnote}{\fnsymbol{footnote}}
\footnotetext[1]{Corresponding author}
\endgroup
\vspace{6pt}

\begin{center}
\begin{minipage}{0.90\textwidth}
\footnotesize
\hrule height 0.4pt
\vspace{5pt}
\noindent\textbf{Abstract}\quad
We determine the Morse index and normal nullity of all linearly full
homogeneous minimal two-spheres in complex hyperquadrics.
These invariants satisfy a common relation involving only the ambient
dimension and the induced Gauss curvature.
We also prove that every normal Jacobi field is integrable through minimal immersions.
\par\smallskip
\noindent\textbf{Keywords}\quad
Morse index, normal nullity, homogeneous minimal two-sphere,
complex hyperquadric, Jacobi fields
\par\smallskip
\noindent\textbf{MSC(2020)}\quad 53C42, 53C35, 58E12
\par\vspace{6pt}
\hrule height 0.4pt
\end{minipage}
\end{center}

\section{Introduction}\label{sec:intro}
The study of homogeneous minimal two-spheres is a basic step toward
classifying constant-curvature minimal two-spheres in symmetric spaces
\cite{XuSurvey}. In round spheres and complex projective spaces, every
linearly full constant-curvature minimal two-sphere is homogeneous:
the classical models are the Bor\r{u}vka spheres and
the Veronese sequence, respectively
\cite{Boruvka,Calabi,doCarmoWallach,BandoOhnita,Bolton}.

For previous results on minimal two-spheres in complex hyperquadrics, see
\cite{WolfsonQuadric,BahyWood,WangJiaoTotallyReal,JiaoLiParallel,MaciaNagatomoTakahashi,MaciaNagatomo,NagatomoGrassmann,FeiWangXu}.
These immersions are conformal harmonic maps and arise as classical solutions of
two-dimensional nonlinear $\sigma$-models \cite{BykovKrivorol}.
The constant-curvature case has been investigated in
\cite{JiaoWangQn,JiaoLiClassification,JiaoLiQ5,JiaoLiLiQ4,ChiXieXu}.
For homogeneous minimal two-spheres, Peng--Wang--Xu \cite{PWX}
obtained the full classification in terms of real orthogonal
representations of $\SU(2)$, giving explicit models and their geometric
invariants.

For these classified models, we consider the second
variation of area and deformations through minimal immersions.
The Morse index measures instability under normal variations, and
the normal nullity is the dimension of the normal Jacobi kernel.
The integrability problem asks whether every normal Jacobi field
can be realized as the normal velocity of a family of minimal immersions.
Ohnita \cite{OhnitaStability} used representation theory to study the
index and nullity of totally geodesic minimal submanifolds in compact
symmetric spaces.
Ball--Madnick \cite{BallMadnick} computed the index and nullity of
constant-curvature minimal two-spheres in round spheres.
In the harmonic-map setting, Lemaire--Wood \cite{LemaireWood} proved
that Jacobi fields along harmonic two-spheres in $\CP^2$ are integrable.
For constant-curvature minimal two-spheres in $\CP^n$,
Cui--Huang \cite{CuiHuangCP} computed the index and normal nullity
and proved the integrability of all normal Jacobi fields.
The corresponding questions for homogeneous minimal two-spheres in complex
hyperquadrics require a more involved analysis, since the classification
includes several distinct families.

In this paper, we determine the Morse index and normal nullity of all
linearly full homogeneous minimal two-spheres in complex hyperquadrics
and prove that every normal Jacobi field is integrable through minimal
immersions. The index and normal nullity satisfy a common relation
involving only the ambient dimension and the induced Gauss curvature.
We reduce the normal Jacobi operator to real symmetric matrices and compute their
inertia by congruence transformations.
We prove integrability by deforming the frames and comparing the
dimensions of the parameter spaces with the normal nullity.
We equip $Q_n$ with the metric induced by the Fubini--Study metric
of holomorphic sectional curvature $4$, and write $\Ind f$ and $\Nul f$
for the Morse index and normal nullity, respectively.
The four families in the classification of Peng--Wang--Xu \cite{PWX}
and their parameters are recalled in Theorem~\ref{hom:classification}.

\begin{theorem}[Index and nullity]\label{main:theorem}
Let $f:S^2\to Q_n$ be a linearly full homogeneous minimal immersion
with induced Gauss curvature $K$. Then
\begin{equation}\label{eq:main-unified}
 \Ind f+\frac12\Nul f=n\left(1+\frac4K\right)-3.
\end{equation}
The individual values of $\Ind f$ and $\Nul f$ are given in
Table~\ref{models:formulas}.
\end{theorem}

\begin{table}[H]
\centering
\caption{Index and normal nullity of the classified families.
The parameters are as in Theorem~\ref{hom:classification}, with $\lambda>0$
in Type~I and $\mu>0$ in Type~II.
Here $n$ is the dimension of the ambient quadric in the corresponding row.}
\label{models:formulas}
\renewcommand{\arraystretch}{1.35}
\begin{tabularx}{\linewidth}{@{}p{.12\linewidth}p{.24\linewidth}X@{}}
\toprule
Type & Ambient quadric & Index and normal nullity \\
\midrule
I & $Q_{2n_0}$ &
$\begin{aligned}
 \Ind f&=(n_0-\lambda)\{n_0(n_0+\lambda)-1\},\\
 \Nul f&=4n_0^2+6n_0-2\lambda-6.
\end{aligned}$ \\
\midrule
II & $\begin{aligned}
 &Q_{m_0-1}\quad(\tau=0),\\
 &Q_{2m_0}\quad(\tau>0)
\end{aligned}$ &
$\begin{aligned}
 \Ind f&=\tfrac{m_0-\mu}{2}\{n(m_0+\mu)-2\},\\
 \Nul f&=2n(m_0+1)+2(m_0-\mu)-6.
\end{aligned}$ \\
\midrule
III.I & $Q_{2m}$ &
$\begin{aligned}
 \Ind f&=m(m^2-1),\\
 \Nul f&=4m^2+6m-6.
\end{aligned}$ \\
\midrule
III.II & $Q_{m_1+m_2}$ &
$\begin{aligned}
 \Ind f&=\tfrac{n-2}{4}(m_1^2+m_2^2+n),\\
 \Nul f&=m_1^2+m_2^2+\tfrac{n(n+6)}2-6.
\end{aligned}$ \\
\bottomrule
\end{tabularx}
\end{table}

The formulas in Table~\ref{models:formulas} show that the stable models
are precisely Type~I with $\lambda=n_0$, Type~II with $\mu=m_0$,
and Type~III.II with $(m_1,m_2)=(2,0)$.
The last model is totally real and has $K=2$.

Our second result realizes all normal Jacobi fields as normal velocities
in a smooth family of minimal immersions.

\begin{theorem}[Integrability of normal Jacobi fields]\label{main:integrability}
Let $f:S^2\to Q_n$ be a linearly full homogeneous minimal immersion,
and let $J_f$ be its normal Jacobi operator.
There exist a neighbourhood $U$ of $0$ in $\ker J_f$ and a smooth map
$\Phi:U\times S^2\to Q_n$ such that each $\Phi(v,\cdot)$ is a conformal
minimal immersion with respect to the complex structure on $S^2$ induced by $f$,
and
\[
  \Phi(0,\cdot)=f,\qquad
  \left(\left.\frac{d}{dt}\right|_{0}\Phi(tv,\cdot)\right)^\perp=v
       \quad(v\in\ker J_f).
\]
\end{theorem}

The paper is organized as follows.
Section~\ref{sec:setup} recalls the classification and moving-frame
formulas for homogeneous minimal two-spheres in $Q_n$.
Section~\ref{sec:index} gives the mode decomposition and explicit
matrix blocks of the normal Jacobi operator, and derives the index and nullity formulas in
Theorem~\ref{main:theorem}.
Section~\ref{sec:integrability} constructs families of minimal immersions
and proves the integrability result in Theorem~\ref{main:integrability}.
Appendix~\ref{app:congruence} proves Proposition~\ref{index:inertia}.

\section{Homogeneous minimal two-spheres}\label{sec:setup}
We write $Q_n=\{[z]\in\CP^{n+1}:zz^t=0\}$.
An immersion $f:S^2\to Q_n$ is homogeneous if a subgroup of the
ambient isometry group acts transitively on its image.
It is linearly full if no real orthogonal transformation maps its image
into a standard subquadric $Q_m\subset Q_n$ with $m<n$.
We fix the representation and moving-frame notation and recall the
homogeneous models of Peng--Wang--Xu \cite{PWX}.

The special unitary group is
\[
 \SU(2)=\left\{\begin{pmatrix}a&b\\-\bar b&\bar a\end{pmatrix}:
 |a|^2+|b|^2=1\right\}.
\]
Its Lie algebra is
\[
 \mathfrak{su}(2)=\left\{\begin{pmatrix}\mathrm{i}x&y\\-\bar y&-\mathrm{i}x\end{pmatrix}:
 x\in\R,\ y\in\C\right\}.
\]
The basis
\[
 \varepsilon_1=\begin{pmatrix}\mathrm{i}&0\\0&-\mathrm{i}\end{pmatrix},\quad
 \varepsilon_2=\begin{pmatrix}0&1\\-1&0\end{pmatrix},\quad
 \varepsilon_3=\begin{pmatrix}0&\mathrm{i}\\\mathrm{i}&0\end{pmatrix}
\]
satisfies
$[\varepsilon_1,\varepsilon_2]=2\varepsilon_3$,
$[\varepsilon_2,\varepsilon_3]=2\varepsilon_1$ and
$[\varepsilon_3,\varepsilon_1]=2\varepsilon_2$.
The right Maurer--Cartan form and its structure equation are
\begin{equation}\label{rep:MC-structure}
 \Theta=\dd g\,g^{-1}
 =\begin{pmatrix}\mathrm{i}\omega&\varphi\\-\bar\varphi&-\mathrm{i}\omega\end{pmatrix},
 \qquad \dd\Theta=\Theta\wedge\Theta,\qquad
 \dd(\mathrm{i}\omega)=-\varphi\wedge\bar\varphi,\quad
 \dd\varphi=2\mathrm{i}\omega\wedge\varphi,
\end{equation}
where $\omega$ is real and $\varphi$ is complex.

Let $V_n$ be the $(n+1)$-dimensional complex vector space of homogeneous
polynomials of degree $n$ in $z_0,z_1$, where $n\ge0$.
For $f_1=\sum_{k=0}^n a_kz_0^kz_1^{n-k}$ and
$f_2=\sum_{k=0}^n b_kz_0^kz_1^{n-k}$, we define
\[
 \langle f_1,f_2\rangle_\C=\sum_{k=0}^n a_k\bar b_k\,k!(n-k)!.
\]
Then $V_n$ has the orthonormal complex basis
\[
 v_{k,n}=\frac{z_0^kz_1^{n-k}}{\sqrt{k!(n-k)!}},
 \qquad 0\le k\le n.
\]
In row-vector notation, we define the unitary representation $\rho_n$ by
\begin{equation}\label{common:right-representation}
 (f\rho_n(g))(z_0,z_1)=f((z_0,z_1)g^{-1})
 =f(\bar a z_0+\bar b z_1,-bz_0+az_1).
\end{equation}
The differential representation is defined by
\[
 v_{k,n}\,\dd\rho_n(\varepsilon)
 =\left.\frac{\dd}{\dd t}\bigl(v_{k,n}\rho_n(\exp(t\varepsilon))\bigr)\right|_{t=0},
 \qquad \varepsilon\in\mathfrak{su}(2).
\]
We set $a_{k,n}=\sqrt{(k+1)(n-k)}$. Differentiation gives
\begin{equation}\label{rep:ladder}
 \begin{gathered}
 v_{k,n}\dd\rho_n(\varepsilon_1)=\mathrm{i}(n-2k)v_{k,n},\qquad
 v_{k,n}\dd\rho_n(\varepsilon_2)=a_{k-1,n}v_{k-1,n}-a_{k,n}v_{k+1,n},\\
 v_{k,n}\dd\rho_n(\varepsilon_3)=-\mathrm{i}(a_{k-1,n}v_{k-1,n}+a_{k,n}v_{k+1,n}).
 \end{gathered}
\end{equation}
Terms with indices outside $0\le k\le n$ are zero.
For $0\le k\le n$, we set
\[
 \lambda_k=n-2k,\qquad v_{\lambda_k,n}:=v_{k,n},\qquad
 a_{\lambda_k,n}:=a_{k,n}
 =\tfrac12\sqrt{(n+\lambda_k)(n-\lambda_k+2)}.
\]
We use the weight notation from this point on.
The set of weights is $\Delta_n=\{\lambda_0,\lambda_1,\ldots,\lambda_n\}$,
with highest weight $\lambda_0=n$.
We set $a_{\nu,n}=0$ for $\nu\notin\Delta_n$; in particular,
$a_{-n,n}=a_{n+2,n}=0$.
The weight-space decomposition is
\[
 V_n=\bigoplus_{k=0}^n V_{\lambda_k,n},\qquad
 V_{\lambda_k,n}=\C v_{\lambda_k,n}.
\]
The weight spaces are mutually orthogonal.
The representations $(V_n,\rho_n)$, $n\ge0$, are pairwise
inequivalent and exhaust the irreducible unitary representations
of $\SU(2)$.  Every finite-dimensional unitary representation is
an orthogonal sum of these representations.

To describe the real representations, we define the conjugate-linear map
$\mathcal J$ on $V_n$ by
\[
 (\mathcal J f)(z_0,z_1)=\overline{f(-\bar z_1,\bar z_0)}.
\]
It satisfies $\mathcal Jv_{\lambda_k,n}=(-1)^kv_{-\lambda_k,n}$ and
\begin{equation}\label{rep:structure-map}
 \mathcal J(v\rho_n(g))=(\mathcal Jv)\rho_n(g),\qquad
 \mathcal J^2=(-1)^nI.
\end{equation}
Restricting the scalar field from $\C$ to $\R$ gives the real
vector space $(V_n)_\R$, of dimension $2(n+1)$, with inner product
$\langle v,w\rangle=\operatorname{Re}\langle v,w\rangle_\C$.
Inner products without a subscript denote the corresponding real inner products.
The resulting real representation, denoted by $((V_n)_\R,\rho_n)$,
is orthogonal with respect to this inner product.
We set $u_{\lambda_k,n}=\mathrm{i}v_{\lambda_k,n}$. An orthonormal real basis is
$\{v_{\lambda_k,n},u_{\lambda_k,n}:0\le k\le n\}$, and
\[
 v_{\lambda,n}\dd\rho_n(\varepsilon_1)=\lambda u_{\lambda,n},\qquad
 u_{\lambda,n}\dd\rho_n(\varepsilon_1)=-\lambda v_{\lambda,n}.
\]

For odd $n$, the real representation $(V_n)_\R$ is irreducible.
For even $n$, it is reducible, with the orthogonal decomposition
\[
 (V_n)_\R=(I+\mathcal J)(V_n)_\R
 \oplus(I-\mathcal J)(V_n)_\R.
\]
The restrictions of $\rho_n$ to the two summands are equivalent
irreducible real representations.
For even $n$, we set
\begin{equation}\label{rep:real-basis}
 w_{\lambda_k,n}=\frac{v_{\lambda_k,n}+(-1)^{(n-\lambda_k)/2}u_{-\lambda_k,n}}{\sqrt2},
 \qquad
 \widetilde V_n=\operatorname{span}_\R\{w_{\lambda_k,n}:0\le k\le n\}.
\end{equation}
This real subspace is invariant under $\rho_n$.
We denote the restriction of $\rho_n$ to $\widetilde V_n$ by
$\widetilde\rho_n$. It is equivalent to the restriction of
$\rho_n$ to $(I+\mathcal J)(V_n)_\R$.
The vectors $w_{\lambda_k,n}$ form an orthonormal real basis of
$\widetilde V_n$, so $\dim_\R\widetilde V_n=n+1$.
Its zero-weight line is $\R w_{0,n}$, whereas the odd-degree
real representations have no zero weight.

The representations $((V_{2n'+1})_\R,\rho_{2n'+1})$ and
$(\widetilde V_{2n'},\widetilde\rho_{2n'})$, $n'\ge0$, are pairwise
inequivalent and exhaust the irreducible real representations of
$\SU(2)$.  Their real dimensions are $4n'+4$ and $2n'+1$, respectively.
Thus every finite-dimensional real representation $(V,\rho)$,
equipped with an invariant inner product, has an orthogonal decomposition
\[
 V=\bigoplus_\alpha(V_{n_\alpha})_\R
      \oplus\bigoplus_\beta\widetilde V_{m_\beta},\qquad
 \rho=\bigoplus_\alpha\rho_{n_\alpha}
      \oplus\bigoplus_\beta\widetilde\rho_{m_\beta},
\]
where the $n_\alpha$ are odd and the $m_\beta$ are even.
Any real basis of $V$ is a complex basis of its complexification
$V^\C=V\otimes_\R\C$.  The same real matrices define the complex-linear
extension of $\rho$, again denoted by $\rho$.
In the orbit formulas below, $\mathrm{i}$ denotes scalar multiplication in $V^\C$.

The differential $\dd\rho$ extends complex linearly to
$\mathfrak{sl}(2,\C)$ acting on $V^\C$.
For the standard basis
\[
 \sigma_1=\begin{pmatrix}1&0\\0&-1\end{pmatrix},\qquad
 \sigma_2=\begin{pmatrix}0&1\\0&0\end{pmatrix},\qquad
 \sigma_3=\begin{pmatrix}0&0\\1&0\end{pmatrix},
\]
we write
\begin{equation}\label{rep:HAB}
 H=\dd\rho(\sigma_1)=-\mathrm{i}\dd\rho(\varepsilon_1),\qquad
 A=\dd\rho(\sigma_2)=\frac{\dd\rho(\varepsilon_2)-\mathrm{i}\dd\rho(\varepsilon_3)}2,
 \qquad
 B=\dd\rho(\sigma_3)=-\frac{\dd\rho(\varepsilon_2)+\mathrm{i}\dd\rho(\varepsilon_3)}2.
\end{equation}
For $\lambda\in\Delta_n$, their action on the original complex representation $V_n$ is given by
$v_{\lambda,n}H=\lambda v_{\lambda,n}$,
$v_{\lambda,n}A=-a_{\lambda,n}v_{\lambda-2,n}$ and
$v_{\lambda,n}B=-a_{\lambda+2,n}v_{\lambda+2,n}$.
The commutator relations are
\[
 [H,A]=2A,\qquad [H,B]=-2B,\qquad [A,B]=H.
\]
The matrix-valued Maurer--Cartan form is therefore
\begin{equation}\label{rep:representation-MC}
 \dd\rho\,\rho^{-1}=H(\mathrm{i}\omega)+A\varphi-B\bar\varphi.
\end{equation}

For $g\in\SU(2)$, we write
\[
 \phi_{\lambda,\alpha}:=v_{\lambda,n_\alpha}\rho_{n_\alpha}(g),\qquad
 \phi'_{\lambda,\alpha}:=u_{\lambda,n_\alpha}\rho_{n_\alpha}(g),\qquad
 \widetilde\phi_{\mu,\beta}:=w_{\mu,m_\beta}\widetilde\rho_{m_\beta}(g).
\]
Here $\alpha$ and $\beta$ label the summands, with
$\lambda\in\Delta_{n_\alpha}$ and $\mu\in\Delta_{m_\beta}$.

Let $T=\{\exp(t\varepsilon_1):t\in\R\}$ be a maximal
torus of $\SU(2)$.

By Peng--Wang--Xu \cite{PWX}, every homogeneous immersion
$f:S^2\to Q_n$ has an orbit representation
\begin{equation}\label{hom:orbit-realization}
 f([g])=[z_0\rho(g)],\qquad [g]\in S^2\simeq\SU(2)/T,
\end{equation}
where $\rho:\SU(2)\to\mathrm O(n+2)$ is an orthogonal representation
and $[z_0]\in Q_n$ is fixed by $\rho(T)$.
We normalize $z_0$ so that $z_0\bar z_0^{\,t}=1$.
Since $[z_0]$ is fixed by $\rho(T)$, we may write $z_0H=-\lambda z_0$
for some $\lambda\in\mathbb Z$.

For the unit lift $Z=z_0\rho(g)$, we set
\[
 X=z_0A\rho(g),\qquad Y=-z_0B\rho(g).
\]
The row vectors $X,Y$ represent mutually Hermitian-orthogonal tangent
vectors at $[Z]\in Q_n$.
With $\varpi=\langle\dd Z,Z\rangle_\C=-\mathrm{i}\lambda\omega$,
the moving-frame formulas of \cite{PWX} give
\begin{equation}\label{common:lift-differential}
 \dd Z=\varpi Z+X\varphi+Y\bar\varphi,\qquad
 \dd s^2=(|X|^2+|Y|^2)\varphi\bar\varphi
        =\frac4K\varphi\bar\varphi,
\end{equation}
where $K$ is the Gauss curvature of the induced metric.
The induced Levi--Civita connection satisfies
$\nabla\varphi=2\mathrm{i}\omega\otimes\varphi$.
We define the covariant differentials
\begin{equation}\label{common:tangent-covariant-definition}
 DX=\dd X-\varpi X+2\mathrm{i}\omega X,
 \qquad DY=\dd Y-\varpi Y-2\mathrm{i}\omega Y.
\end{equation}
Differentiation using \eqref{rep:representation-MC} gives
\begin{equation}\label{common:tangent-covariant-calculation}
 DX=z_0(A^2\varphi-AB\bar\varphi)\rho(g),\qquad
 DY=z_0(-BA\varphi+B^2\bar\varphi)\rho(g).
\end{equation}
The covariant derivative in $Q_n$ is the Hermitian projection of the
ambient covariant derivative in $\CP^{n+1}$ onto
$\operatorname{span}_\C\{Z,\bar Z\}^{\perp}$.

\Needspace{10\baselineskip}
\begin{theorem}[Homogeneous minimal two-spheres \cite{PWX}]\label{hom:classification}
Let $f:S^2\to Q_n$ be a linearly full homogeneous minimal immersion.
Then $f$ is congruent to one of the following models, each linearly full
in its stated ambient quadric.
\begin{enumerate}[label=\textup{\Roman*.},leftmargin=3.5em,labelwidth=3em,itemsep=3pt]
\item Let $n_1=n_0\ge1$ be odd and $\lambda\in\Delta_{n_0}$.
For $\rho=\rho_{n_0}$,
\[
 [g]\longmapsto[\phi_{\lambda,1}+\mathrm{i}\phi'_{\lambda,1}]
 \in Q_{2n_0}.
\]
Its Gauss curvature is $K=8/(n_0^2+2n_0-\lambda^2)$.
\item Let $m_0\ge2$ be even and $\mu\in\Delta_{m_0}\setminus\{0\}$.
For $0<\tau\le\pi/4$, we take $m_1=m_2=m_0$ and
$\rho=\widetilde\rho_{m_0}\oplus\widetilde\rho_{m_0}$, and set
\[
 [g]\longmapsto[\cos\tau(\widetilde\phi_{\mu,1}
 +\mathrm{i}(-1)^{(m_0-\mu)/2}\widetilde\phi_{-\mu,1}),
 \mathrm{i}\sin\tau(\widetilde\phi_{\mu,2}
 +\mathrm{i}(-1)^{(m_0-\mu)/2}\widetilde\phi_{-\mu,2})]\in Q_{2m_0}.
\]
At $\tau=0$, we take $m_1=m_0$ and
$\rho=\widetilde\rho_{m_0}$, and use the single-summand model
\[
 [g]\longmapsto[\widetilde\phi_{\mu,1}
 +\mathrm{i}(-1)^{(m_0-\mu)/2}\widetilde\phi_{-\mu,1}]
 \in Q_{m_0-1}.
\]
In both cases, $K=8/(m_0^2+2m_0-\mu^2)$.
\item[\textup{III.I.}] Let $m_1=m_2=m\ge2$ be even.
For $\rho=\widetilde\rho_m\oplus\widetilde\rho_m$,
\[
 [g]\longmapsto[\widetilde\phi_{0,1},\mathrm{i}\widetilde\phi_{0,2}]
 \in Q_{2m},\qquad K=\frac8{m^2+2m}.
\]
\item[\textup{III.II.}] Let $m_1>m_2\ge0$ be even.
For $\rho=\widetilde\rho_{m_1}\oplus\widetilde\rho_{m_2}$,
\[
 [g]\longmapsto[\widetilde\phi_{0,1},\mathrm{i}\widetilde\phi_{0,2}]
 \in Q_{m_1+m_2},\qquad
 K=\frac{16}{m_1^2+m_2^2+2(m_1+m_2)}.
\]
\end{enumerate}
\end{theorem}

\begin{remark}
Each representative in Theorem~\ref{hom:classification} is taken in the
quadric associated with the real span of the real and imaginary parts
of its lift. For the two-summand Type~I
representative of Peng--Wang--Xu \cite{PWX}, we consider the real linear map
\[
 v_{\nu,n_0}\longmapsto
 (\cos\tau\,v_{\nu,n_0},-\sin\tau\,u_{\nu,n_0}),\qquad
 u_{\nu,n_0}\longmapsto
 (\cos\tau\,u_{\nu,n_0},\sin\tau\,v_{\nu,n_0})
 \quad(\nu\in\Delta_{n_0}).
\]
It is equivariant and isometric, and its complexification sends the
single-summand orbit to the two-summand representative.
For the resulting lift $Z$, real irreducibility gives
\[
 \dim_\R\operatorname{span}_\R
 \{\operatorname{Re}Z(g),\operatorname{Im}Z(g):g\in\SU(2)\}
 =2n_0+2.
\]
Hence the linearly full target is $Q_{2n_0}$.
For Type~II, exchanging the two summands and changing the sign of one
allows us to assume $0\le\tau\le\pi/4$. At $\tau=0$, the real and imaginary
parts of the lift span the nonzero irreducible summand of real dimension
$m_0+1$, so the linearly full target is $Q_{m_0-1}$.
\end{remark}

\begin{remark}\label{hom:central-action}
Since the polynomials in $V_n$ are homogeneous of degree $n$,
\eqref{common:right-representation} gives $\rho_n(-I)=(-1)^nI$.
For the models in Theorem~\ref{hom:classification}, this gives
$\rho(-I)=-I$ in Type~I and $\rho(-I)=I$ in Types~II, III.I and III.II.
Thus $-I$ induces the identity isometry on each of these ambient quadrics.
\end{remark}

For Types~I and II in Section~\ref{sec:integrability}, we use the positive parameters from
Table~\ref{models:formulas} and the notation $(n_0,\lambda)$ for both types,
with $(n_0,\lambda)=(m_0,\mu)$ in Type~II.
We choose an orthonormal weight frame $Z_\nu(g)=Z_\nu(1)\rho(g)$
for the model, where $\nu\in\Delta_{n_0}$. The phases are chosen so that
the weight-frame equations have the signs in \eqref{rep:ladder}.
Here $1=I_2$ denotes the identity of $\SU(2)$,
so the $Z_\nu(1)$ are the fixed initial frame vectors.
After precomposing with the antipodal map if
necessary, we choose the orientation for which $[Z_{n_0}]$
is holomorphic and $f=[Z_\lambda]$.
The weight-frame equations are
\begin{equation}\label{int:nonzero-ladder}
 \dd Z_\nu=\mathrm i\nu\omega Z_\nu
 -a_{\nu,n_0}Z_{\nu-2}\varphi
 +a_{\nu+2,n_0}Z_{\nu+2}\bar\varphi.
\end{equation}

For Types~III.I and III.II, we write $n=m_1+m_2$, where $m_1\ge m_2\ge0$ are even and $m_1>0$.
The inclusion $\widetilde V_{m_i}\subset V_{m_i}$ extends complex linearly
to an isomorphism $\widetilde V_{m_i}^{\C}\cong V_{m_i}$.
Under this identification, for $i=1,2$ and
integers $-m_i/2\le a\le m_i/2$, we use the frames
\begin{equation}\label{zero:weight-frames}
 f_{i,a}=\frac{1+\mathrm{i}(-1)^{m_i/2}}{\sqrt2}\,
 v_{2a,m_i}\rho_{m_i}(g),\qquad
 Z=\frac{f_{1,0}+\mathrm{i}f_{2,0}}{\sqrt2}.
\end{equation}
Then $\bar f_{i,a}=(-1)^af_{i,-a}$ and $\varpi=0$.

\section{Morse index and nullity}\label{sec:index}
We compute the Morse index and normal nullity of the models in
Theorem~\ref{hom:classification} by reducing the normal Jacobi operator
to real symmetric matrices.
Section~\ref{sec:jacobi-operator} derives the index form and expresses
the index and normal nullity in terms of the eigenvalue counts of these matrices.
Section~\ref{sec:jacobi-matrices} computes the matrices for each family
and proves the formulas in Theorem~\ref{main:theorem}.

\subsection{Jacobi operator and mode decomposition}\label{sec:jacobi-operator}

Let $f:S^2\to Q_n$ be a linearly full homogeneous minimal immersion
written as in \eqref{hom:orbit-realization}. We identify $TS^2$ with its image
under $\dd f$ and let $e_1,e_2$ be a local orthonormal tangent frame
for the induced metric.

For normal fields $V,W$, the Jacobi operator and the second variation of
area are given by
\begin{equation}\label{common:Jacobi-convention}
 J=\Delta^\perp+\mathcal R+\mathcal B,\qquad
 Q_f(V,W)=-\int_{S^2}\langle JV,W\rangle\,\dd A,
\end{equation}
where
\begin{equation}\label{common:normal-Laplacian}
 \Delta^\perp V=\sum_{a=1}^2\left(
 \nabla^\perp_{e_a}\nabla^\perp_{e_a}V
 -\nabla^\perp_{\nabla_{e_a}e_a}V\right),
\end{equation}
\begin{equation}\label{common:Jacobi-geometric-terms}
 \mathcal R(V)=\sum_{a=1}^2\bigl(R(V,e_a)e_a\bigr)^\perp,\qquad
 \mathcal B(V)=\sum_{a,b=1}^2
       \langle V,\mathrm{II}_{ab}\rangle\mathrm{II}_{ab}.
\end{equation}
Here $\nabla$ is the Levi--Civita connection of the induced metric,
$\nabla^\perp$ is the normal connection, $R$ is the curvature tensor of
$Q_n$, and $\mathrm{II}_{ab}=\mathrm{II}(e_a,e_b)$.
Thus $\Ind f$ is the number of positive eigenvalues of $J$, counted with
their real multiplicities, and $\Nul f=\dim_\R\ker J$ \cite{Simons,Urbano}.

Relative to the unit lift $Z$, we represent tangent and normal vectors
by complex row vectors, with real inner product
$\langle V,W\rangle=\operatorname{Re}(V\bar W^t)$.
In this representation, the normal connection is given by
\[
 \nabla^\perp V=(\dd V-\varpi V)^\perp,
\]
where $(\cdot)^\perp$ denotes the orthogonal projection onto the normal
space of $f$.
Using \eqref{common:lift-differential} and
\eqref{common:tangent-covariant-calculation}, we choose the local
orthonormal tangent frame $e_1,e_2$ as follows:
\begin{equation}\label{common:tangent-frame}
 e_1=\frac{\sqrt K}{2}z_0(A-B)\rho(g),\qquad
 e_2=\frac{\mathrm{i}\sqrt K}{2}z_0(A+B)\rho(g),\qquad
 \sum_{a=1}^2 e_ae_a^t=\frac K2z_0(AB+BA)z_0^t.
\end{equation}
For a normal field $V$, the second fundamental form satisfies
\[
 \langle V,\mathrm{II}_{11}\rangle
 =\frac K4\langle V,z_0(A^2+B^2)\rho(g)\rangle,\qquad
 \langle V,\mathrm{II}_{12}\rangle
 =\frac K4\langle V,\mathrm{i}z_0(A^2-B^2)\rho(g)\rangle.
\]
For the curvature calculation, extend $Z$ to a local unit lift on $Q_n$
and $V,e_a$ to tangent vector fields on $Q_n$, with $\varpi=\dd Z\bar Z^t$.
The identities
\[
 ZZ^t=0,\qquad Z\bar Z^t=1,\qquad VZ^t=V\bar Z^t=0
\]
give, by differentiation and orthogonal projection, the covariant
differential of $V$ in $Q_n$:
\[
 \dd V-\varpi V+(V\dd\bar Z^t)Z+(V\dd Z^t)\bar Z.
\]
Since $\dd\varpi=-\dd Z\wedge\dd\bar Z^t$ and
$\operatorname{Re}(V\bar e_a^t)=0$ along $f$, we have
\[
 -\dd\varpi(V,e_a)=V\bar e_a^t-e_a\bar V^t=2V\bar e_a^t
 \qquad(a=1,2).
\]
Applying the definition of curvature to this covariant differential,
and using $|e_a|=1$ along $f$, gives
\[
 R(V,e_a)e_a
 =V-\bigl[\dd\varpi(V,e_a)+e_a\bar V^t\bigr]e_a
   +(e_ae_a^t)\bar V-(Ve_a^t)\bar e_a
 =V+3(V\bar e_a^t)e_a+(e_ae_a^t)\bar V-(Ve_a^t)\bar e_a.
\]
Since $|V\bar e_a^t|^2=\langle V,\mathrm i e_a\rangle^2$,
taking real inner products with $V$ and summing over $a=1,2$ gives
the first identity below. The second follows from the orthogonal decomposition
of $\dd V-\varpi V$ along $f$:
\[
\begin{aligned}
 \sum_{a=1}^2\langle R(V,e_a)e_a,V\rangle
 &={}2|V|^2+3\sum_{a=1}^2\langle V,\mathrm{i}e_a\rangle^2
      -\sum_{a=1}^2|Ve_a^t|^2
      +\operatorname{Re}\left[
       \left(\sum_{a=1}^2 e_ae_a^t\right)\overline{VV^t}\right],\\
 |\dd V-\varpi V|^2
 &={} |\nabla^\perp V|^2
      +\sum_{a,b=1}^2\langle V,\mathrm{II}_{ab}\rangle^2
      +\sum_{a=1}^2\left(\langle V,\mathrm{i}e_a\rangle^2
                         +|Ve_a^t|^2\right).
\end{aligned}
\]
With $\mathrm{II}_{22}=-\mathrm{II}_{11}$, substitution in
\eqref{common:Jacobi-convention}--\eqref{common:Jacobi-geometric-terms} yields
\begin{equation}\label{common:index-form}
 Q_f(V,V)=\int_{S^2}\Bigl\{|\dd V-\varpi V|^2-2|V|^2
 -4\sum_{a=1}^2\bigl(\langle V,\mathrm{i}e_a\rangle^2
                         +\langle V,\mathrm{II}_{1a}\rangle^2\bigr)
 -\operatorname{Re}\left[\left(\sum_{a=1}^2 e_ae_a^t\right)
                    \overline{VV^t}\right]\Bigr\}\,\dd A.
\end{equation}
Norms of one-forms are computed with respect to the induced metric.

The right $\SU(2)$ action is induced by ambient isometries and preserves
the normal connection, curvature and second fundamental form.
It therefore commutes with $J$.
By Remark~\ref{hom:central-action}, $-I$ acts as the identity on $L^2(Nf)$.
Since it acts as $(-1)^mI$ on $V_m$, only even highest weights occur,
with real irreducible representations $\widetilde V_{2\ell}$.

Since $\dim_\R\widetilde V_{2\ell}=2\ell+1$ is odd, every real
equivariant endomorphism $L$ of $\widetilde V_{2\ell}$ has a real
eigenvalue $c$. The nonzero invariant subspace $\ker(L-cI)$ is the
whole space by irreducibility. Thus
\[
 \operatorname{End}_{\SU(2)}(\widetilde V_{2\ell})=\R I.
\]
Let $N_Tf\subset T_{[z_0]}Q_n$ be the normal space of $f$ at
the base coset $T\in\SU(2)/T$, where $f(T)=[z_0]$.
The isotropy group $T$ acts linearly on $N_Tf$.
All $\Hom$ spaces below consist of real linear equivariant maps.
Evaluation at $T$ gives the Frobenius reciprocity isomorphism
\cite{BrockerTomDieck,Helgason}:
\[
 \Hom_{\SU(2)}\bigl(\widetilde V_{2\ell},C^\infty(Nf)\bigr)
 \xrightarrow{\;\cong\;}\Hom_T(\widetilde V_{2\ell},N_Tf),\qquad
 \bigl[u\mapsto\Phi(u)\bigr]
 \longmapsto\bigl[u\mapsto\Phi(u)(T)\bigr].
\]
This isomorphism and the Peter--Weyl theorem for the real
Hilbert space $L^2(Nf)$ give the orthogonal decomposition
\cite{BrockerTomDieck}
\begin{equation}\label{common:normal-Peter-Weyl}
 L^2(Nf)\cong\widehat{\bigoplus}_{\ell\ge0}
 \widetilde V_{2\ell}\otimes_\R\Hom_T(\widetilde V_{2\ell},N_Tf).
\end{equation}
Since $\SU(2)$ acts only on the first factor and
$\operatorname{End}_{\SU(2)}(\widetilde V_{2\ell})=\R I$,
equivariance and self-adjointness of $J$ give the orthogonal decomposition
\begin{equation}\label{common:Jacobi-decomposition}
 J\cong\bigoplus_{\ell\ge0}\bigl(I_{2\ell+1}\otimes J_\ell\bigr),
\end{equation}
where $J_\ell$ is a real symmetric matrix in orthonormal coordinates on
$\Hom_T(\widetilde V_{2\ell},N_Tf)$.
Thus each eigenvalue of $J_\ell$ with multiplicity $m$ occurs with
multiplicity $(2\ell+1)m$ in the corresponding summand.

To compute the matrices $J_\ell$, we use the following expansion of normal fields.
By fixed complex linear combinations of the real frames in
Section~\ref{sec:setup}, we obtain a complex orthonormal weight frame
$E_1,\ldots,E_{n+2}$; the explicit choices are given in
Section~\ref{sec:jacobi-matrices}.
For a normal field $V$, write
\[
 V(g)=\sum_{a=1}^{n+2}\eta_a(g)E_a(g),\qquad
 \eta_a:\SU(2)\to\C.
\]
If the weight of $E_a$ minus the weight of $Z$ is $2r_a$,
then its coefficient satisfies
\[
 \eta_a(\exp(t\varepsilon_1)g)=e^{-2\mathrm{i}r_at}\eta_a(g),\qquad r_a\in\mathbb Z.
\]
For $\eta=\eta_a$ and $r=r_a$, the scalar derivatives $\eta A,\eta B$ are given by
\begin{equation}\label{common:coefficient-HAB}
 \dd\eta+2\mathrm{i}r\omega\eta=(\eta A)\varphi-(\eta B)\bar\varphi.
\end{equation}
For fixed $r$, the product of one function satisfying this transformation
law with the conjugate of another descends to $S^2$.
Integration against the induced area form defines the $L^2$ inner product.
An orthonormal basis of this $L^2$ space is
\[
 \{Y_{\ell,j}^{(r)}:\ell\ge|r|,\quad 1\le j\le2\ell+1\},
\]
where $Y_{\ell,j}^{(r)}$ is a normalized matrix coefficient of
$\rho_{2\ell}$. Every such $\eta$ has an expansion with constant coefficients,
\[
 \eta=\sum_{\ell\ge|r|}\sum_{j=1}^{2\ell+1}
 \eta_{\ell,j}Y_{\ell,j}^{(r)},\qquad \eta_{\ell,j}\in\C.
\]
The index $j$ labels the vectors of the real basis in
\eqref{rep:real-basis}. The phases are chosen compatibly with
\eqref{rep:ladder} and \eqref{rep:structure-map}, so that
$\overline{Y_{\ell,j}^{(r)}}=(-1)^{r-1}Y_{\ell,j}^{(-r)}$ and
\begin{equation}\label{common:scalar-mode-ladders}
 Y_{\ell,j}^{(r)}A=a_{2r+2,2\ell}Y_{\ell,j}^{(r+1)},\qquad
 Y_{\ell,j}^{(r)}B=a_{2r,2\ell}Y_{\ell,j}^{(r-1)}.
\end{equation}

For fixed $\ell\ge0$ and $1\le j\le2\ell+1$, write
\[
 \eta_a=(x_a+\mathrm i y_a)Y_{\ell,j}^{(r_a)},\qquad
 x_a,y_a\in\R,\qquad |r_a|\le\ell.
\]
The corresponding $\C^{n+2}$-valued field is
\[
 V_{\ell,j}[x,y]
 =\sum_{|r_a|\le\ell}(x_a+\mathrm i y_a)Y_{\ell,j}^{(r_a)}E_a
 =\sum_{|r_a|\le\ell}\bigl\{
 x_a\bigl(Y_{\ell,j}^{(r_a)}E_a\bigr)
 +y_a\bigl(\mathrm iY_{\ell,j}^{(r_a)}E_a\bigr)\bigr\}.
\]
Thus $x,y$ are real coefficient vectors. Orthonormality of the frame
and scalar modes gives
\[
 \int_{S^2}|V_{\ell,j}[x,y]|^2\,\dd A
 =\sum_{|r_a|\le\ell}(x_a^2+y_a^2)=x^tx+y^ty.
\]
Substituting this expansion into the normality conditions
\[
 V\bar Z^t=VZ^t=0,\qquad
 \operatorname{Re}(V\bar e_1^t)=\operatorname{Re}(V\bar e_2^t)=0
\]
and using the conjugation rule for $Y_{\ell,j}^{(r)}$ gives real linear
constraints on $x,y$, independent of $j$.
These constraints define a real subspace of the $(x,y)$-space,
with its Euclidean inner product. Choose an orthonormal basis of this
subspace and let $\mathbf c$ be the coordinate vector.
Denote the corresponding normal field by $V_{\ell,j}[\mathbf c]$. Then
\[
 \int_{S^2}|V_{\ell,j}[\mathbf c]|^2\,\dd A=\mathbf c^t\mathbf c.
\]
Since $j$ indexes the real basis in \eqref{rep:real-basis}, the
matrix-coefficient transformation law gives
\[
 V_{\ell,j}[\mathbf c]\cdot h
 =\sum_{k=1}^{2\ell+1}
 \bigl(\widetilde\rho_{2\ell}(h)\bigr)_{jk}
 V_{\ell,k}[\mathbf c],\qquad h\in\SU(2),
\]
where $\cdot h$ denotes the induced right action on normal fields.
Together with completeness of the scalar modes, this realizes
\eqref{common:normal-Peter-Weyl}: $j$ labels a basis of the first factor,
while $\mathbf c$ gives orthonormal coordinates on the second.
By \eqref{common:Jacobi-decomposition}, $J$ acts on $\mathbf c$ through
the same real symmetric matrix $J_\ell$ for every $j$.
Every $V\in L^2(Nf)$ has a unique $L^2$-convergent expansion, and
Parseval's identity gives
\begin{equation}\label{common:normal-mode-expansion}
 V=\sum_{\ell\ge0}\sum_{j=1}^{2\ell+1}V_{\ell,j}[\mathbf c_{\ell,j}],\qquad
 \int_{S^2}|V|^2\,\dd A=\sum_{\ell\ge0}\sum_{j=1}^{2\ell+1}\mathbf c_{\ell,j}^t\mathbf c_{\ell,j}.
\end{equation}
We write $\nu_+(M)$, $\nu_0(M)$ and $\nu_-(M)$ for the positive, zero and
negative eigenvalue counts of a real symmetric matrix $M$.
Equations~\eqref{common:Jacobi-convention},
\eqref{common:Jacobi-decomposition} and
\eqref{common:normal-mode-expansion} give the following reduction.

\begin{proposition}[Real matrix reduction]\label{common:Jacobi-reduction}
For each model in Theorem~\ref{hom:classification} and each $\ell\ge0$,
the real symmetric matrix $J_\ell$ defined above is independent of $j$
and satisfies
\begin{equation}\label{common:Jacobi-matrix}
 J V_{\ell,j}[\mathbf c]=V_{\ell,j}[J_\ell\mathbf c].
\end{equation}
For real coordinate vectors $\mathbf c,\mathbf d$,
\begin{equation}\label{common:Jacobi-pullback}
 Q_f\bigl(V_{\ell,j}[\mathbf c],V_{\ell,j}[\mathbf d]\bigr)
 =-\mathbf c^tJ_\ell\mathbf d.
\end{equation}
Consequently,
\begin{equation}\label{common:real-multiplicity-counts}
 \Ind f=\sum_{\ell\ge0}(2\ell+1)\nu_+(J_\ell),\qquad
 \Nul f=\sum_{\ell\ge0}(2\ell+1)\nu_0(J_\ell).
\end{equation}
\end{proposition}

\subsection{Jacobi matrices, index and nullity}\label{sec:jacobi-matrices}
We compute the matrices $J_\ell$ by substituting the moving-frame
expansions of normal fields into the index form \eqref{common:index-form}.
Fix  $(\ell,j_0)$, with $1\le j_0\le2\ell+1$, and set $L=\ell(\ell+1)$.

\medskip\noindent\textit{Types I and II.}
Let $(n_0,\lambda)$ be the positive parameters in
Table~\ref{models:formulas}, with $(n_0,\lambda)=(m_0,\mu)$
in Type~II. Set $c=2$, except in Type~II at $\tau=0$,
where $c=1$; then $n=c(n_0+1)-2$.
We construct the weight frame $E_q^{(i)}$ from the real frames in
Section~\ref{sec:setup}, where $i=1,\ldots,c$ labels the complex weight chains.
Set
\[
 \nu=\lambda+2q,\qquad
 -\frac{n_0+\lambda}{2}\le q\le\frac{n_0-\lambda}{2}.
\]
All combinations below are taken in the complexification of the real
ambient representation, with each summand viewed as a subspace of the
ambient space. In Type~I, define
\[
 E_q^{(1)}=\frac{\phi_{\nu,1}+\mathrm i\phi'_{\nu,1}}{\sqrt2},\qquad
 E_q^{(2)}=\frac{(-1)^{q+\lambda}}{\sqrt2}
       \bigl(\phi_{-\nu,1}-\mathrm i\phi'_{-\nu,1}\bigr).
\]
Thus the index $i$ labels the two complex weight chains arising from
one real irreducible summand. In Type~II with $\tau>0$, define
\[
 \begin{pmatrix}E_q^{(1)}\\[2pt]E_q^{(2)}\end{pmatrix}
 =\frac1{\sqrt2}
 \begin{pmatrix}
  \cos\tau&\mathrm i\sin\tau\\
  \sin\tau&-\mathrm i\cos\tau
 \end{pmatrix}
 \begin{pmatrix}
  \widetilde\phi_{\nu,1}
    +\mathrm i(-1)^{(n_0-\nu)/2}\widetilde\phi_{-\nu,1}\\[2pt]
  \widetilde\phi_{\nu,2}
    +\mathrm i(-1)^{(n_0-\nu)/2}\widetilde\phi_{-\nu,2}
 \end{pmatrix}.
\]
At $\tau=0$, only the first real summand is present, and we set
\[
 E_q^{(1)}=\frac{\widetilde\phi_{\nu,1}
       +\mathrm i(-1)^{(n_0-\nu)/2}\widetilde\phi_{-\nu,1}}{\sqrt2}.
\]
These vectors form an orthonormal weight frame and satisfy
$E_q^{(i)}(g)=E_q^{(i)}(1)\rho(g)$. By \eqref{rep:ladder},
\[
 \dd E_q^{(i)}=-\mathrm i\nu\omega E_q^{(i)}
       +a_{\nu+2,n_0}E_{q+1}^{(i)}\varphi
       -a_{\nu,n_0}E_{q-1}^{(i)}\bar\varphi.
\]
In both types, $E_0^{(1)}=Z$. In Type~I,
$\bar Z=E_{-\lambda}^{(2)}$, whereas in Type~II with $\tau>0$,
\[
 \bar Z=-\mathrm i(-1)^{(n_0-\lambda)/2}
 \bigl(\cos(2\tau)E_{-\lambda}^{(1)}
          +\sin(2\tau)E_{-\lambda}^{(2)}\bigr).
\]
At $\tau=0$, the latter relation reduces to
$\bar Z=-\mathrm i(-1)^{(n_0-\lambda)/2}E_{-\lambda}^{(1)}$.

The weight of $E_q^{(i)}$ relative to $Z$ is $-2q$, so its coefficient
in the fixed mode $(\ell,j_0)$ is a constant multiple of
$Y_{\ell,j_0}^{(-q)}$. This scalar mode exists only for $|q|\le\ell$.
The expansion is therefore
\[
 V=\sum_{i=1}^c\sum_{q=-p}^{r}\eta_q^{(i)}E_q^{(i)},\qquad
 p=\min\{(n_0+\lambda)/2,\ell\},\qquad
 r=\min\{(n_0-\lambda)/2,\ell\}.
\]
Throughout this subsection, coordinates outside the indicated ranges
are zero, and terms involving absent coordinates are omitted.
The condition $V\bar Z^t=0$ gives
$\eta_0^{(1)}=0$. We choose the signs in the following real coordinates
so that the constraints on $x$ and $y$ take the same form:
\begin{equation}\label{matrix:normal-mode}
 \eta_q^{(1)}=(-1)^qY_{\ell,j_0}^{(-q)}
 \begin{cases}
 x_q^{(1)}-\mathrm{i}y_q^{(1)},&q<0,\\
 x_q^{(1)}+\mathrm{i}y_q^{(1)},&q>0,
 \end{cases}
 \qquad
 \eta_q^{(2)}=(-1)^q(x_q^{(2)}-\mathrm{i}y_q^{(2)})Y_{\ell,j_0}^{(-q)}.
\end{equation}
The real constants $x_q^{(i)},y_q^{(i)}$ form coefficient vectors $x,y$
with blocks $x^{(i)},y^{(i)}$.
The conditions $V\bar Z^t=0$ and $\langle V,e_a\rangle=0$ ($a=1,2$)
are equivalent to the following equations, imposed separately on
$u=x$ and $u=y$:
\[
 u_0^{(1)}=0,\qquad
 a_{\lambda,n_0}u_{-1}^{(1)}-a_{\lambda+2,n_0}u_1^{(1)}=0.
\]
For $\ell\ge\lambda$, the expressions for $\bar Z$ show that the remaining
normality condition $VZ^t=0$ is equivalent to
\[
 u_{-\lambda}^{(2)}=0\quad\text{in Type I},\qquad
 \cos(2\tau)u_{-\lambda}^{(1)}+\sin(2\tau)u_{-\lambda}^{(2)}=0
 \quad\text{in Type II}.
\]
This constraint is independent of the preceding two: it involves the
second chain in Type~I and an index $q=-\lambda\le-2$ in Type~II.
For $\ell<\lambda$, the coordinates with $q=-\lambda$ are absent,
so $VZ^t=0$ holds automatically.
The vectors $x$ and $y$ therefore belong to the same real coefficient space.

Since $z_0(AB+BA)=4z_0/K$ and $z_0z_0^t=0$,
\eqref{common:tangent-frame} gives
\[
 \sum_{a=1}^2e_ae_a^t
 =\frac K2z_0(AB+BA)z_0^t=2z_0z_0^t=0.
\]
Here and below, $a=1,2$ indexes the orthonormal tangent frame.
At $\ell=0$, removing the coefficient along $Z$ leaves $2(c-1)$
real parameters. The frame equations give
\[
 \int_{S^2}|\dd V-\varpi V|^2\,\dd A
 =\frac K2\bigl(a_{\lambda,n_0}^2+a_{\lambda+2,n_0}^2\bigr)
   \int_{S^2}|V|^2\,\dd A
 =2\int_{S^2}|V|^2\,\dd A.
\]
Moreover, $\langle V,\mathrm i e_a\rangle=\langle V,\mathrm{II}_{1a}\rangle=0$
for $a=1,2$. Thus \eqref{common:index-form} vanishes on this mode.

For $\ell\ge1$, first take one $x$-block, suppress the superscript $i$,
and write $V=\sum_{q=-p}^r(-1)^qu_qY_{\ell,j_0}^{(-q)}E_q$, where $u_q=x_q$.
The frame and scalar ladder formulas give
\[
 \begin{array}{c@{\qquad}l}
 \text{component}&\text{coefficient in }\dd V-\varpi V\\[3pt]
 E_q\varphi&
 (-1)^q\bigl(a_{2q,2\ell}u_q-a_{\lambda+2q,n_0}u_{q-1}\bigr)
 Y_{\ell,j_0}^{(-q+1)}\\[3pt]
 E_q\bar\varphi&
 (-1)^{q+1}\bigl(a_{2q+2,2\ell}u_q-a_{\lambda+2q+2,n_0}u_{q+1}\bigr)
 Y_{\ell,j_0}^{(-q-1)}.
 \end{array}
\]
Here $q$ ranges over the full weight chain and $u_q=0$ outside $[-p,r]$.
The scalar functions have unit $L^2$ norm and
$|\varphi|^2=|\bar\varphi|^2=K/2$. Thus orthonormality gives
\[
 \begin{aligned}
 \int_{S^2}(|\dd V-\varpi V|^2-2|V|^2)\,\dd A
 ={}&\frac K2\sum_q\Bigl[
 \bigl(a_{2q,2\ell}u_q-a_{\lambda+2q,n_0}u_{q-1}\bigr)^2
 +\bigl(a_{2q+2,2\ell}u_q-a_{\lambda+2q+2,n_0}u_{q+1}\bigr)^2\Bigr]\\
 &\quad-2\sum_{q=-p}^r u_q^2
 =K\left(\sum_{q=-p}^r d_qu_q^2
 -2\sum_{q=-p}^{r-1}b_qu_qu_{q+1}\right)=:Kq_\ell(u),
 \end{aligned}
\]
where
\[
 d_q=\frac12\bigl(a_{2q,2\ell}^2+a_{2q+2,2\ell}^2
 +a_{\lambda+2q,n_0}^2+a_{\lambda+2q+2,n_0}^2\bigr)-\frac2K
 =L-q(2q+\lambda),\qquad 
 b_q=a_{\lambda+2q+2,n_0}a_{2q+2,2\ell}.
\]
The expression for $d_q$ uses $n_0(n_0+2)-\lambda^2=8/K$.
The $y$-blocks give the same form; the sign change in the first chain
in \eqref{matrix:normal-mode} has no effect because $y_0^{(1)}=0$.
Restoring $i$ and summing over the blocks, we obtain for the full field
\[
 \int_{S^2}(|\dd V-\varpi V|^2-2|V|^2)\,\dd A
 =K\sum_{i=1}^c\bigl(q_\ell(x^{(i)})+q_\ell(y^{(i)})\bigr).
\]
The remaining integrals are
\[
 \int_{S^2}\sum_{a=1}^2\langle V,\mathrm i e_a\rangle^2\,\dd A
 =\frac K4\sum_{u\in\{x,y\}}
 \bigl(a_{\lambda,n_0}u_{-1}^{(1)}
       +a_{\lambda+2,n_0}u_1^{(1)}\bigr)^2,
\]
\[
 \int_{S^2}\sum_{a=1}^2\langle V,\mathrm{II}_{1a}\rangle^2\,\dd A
 =\frac{K^2}{16}\sum_{u\in\{x,y\}}
 \bigl(a_{\lambda,n_0}a_{\lambda-2,n_0}u_{-2}^{(1)}
       -a_{\lambda+2,n_0}a_{\lambda+4,n_0}u_2^{(1)}\bigr)^2.
\]
For a coefficient vector $u$ satisfying the normality constraints above,
define
\[
M_\ell(u)
 ={}\sum_{i=1}^cq_\ell(u^{(i)})
 -\bigl(a_{\lambda,n_0}u_{-1}^{(1)}
       +a_{\lambda+2,n_0}u_1^{(1)}\bigr)^2 -\frac K4\bigl(a_{\lambda,n_0}a_{\lambda-2,n_0}u_{-2}^{(1)}
       -a_{\lambda+2,n_0}a_{\lambda+4,n_0}u_2^{(1)}\bigr)^2.
\]

By orthonormality and \eqref{common:index-form}, we obtain
\begin{equation}\label{matrix:nonzero-form}
 \int_{S^2}|V|^2\,\dd A=x^tx+y^ty,\qquad
 Q_f(V,V)=K\{M_\ell(x)+M_\ell(y)\}.
\end{equation}

\medskip\noindent\textit{Types III.I and III.II.}
We use the complex orthonormal weight frames in \eqref{zero:weight-frames}.
Here $n=m_1+m_2$, and $m_1=m_2$ in Type~III.I.
The conditions $V\bar Z^t=VZ^t=0$ remove both zero-weight
components, so there is no $\ell=0$ normal mode. For $\ell\ge1$, we set
$d_i=\min(m_i/2,\ell)$ and $\alpha_i=a_{0,m_i}^2$.
Then $\alpha_1+\alpha_2=4/K$.
The conjugation rules give
\[
 \overline{Y_{\ell,j_0}^{(k)}f_{i,k}}
 =-Y_{\ell,j_0}^{(-k)}f_{i,-k}.
\]
Pairing each term with its conjugate, we write the real and imaginary
parts of $V$ as follows:
\begin{equation}\label{zero:normal-coordinates}
\begin{aligned}
 \operatorname{Re}V&=\sqrt2\operatorname{Re}\left[
 \sum_{k=1}^{d_1}(c_{1,k}+\mathrm{i}d_{1,k})Y_{\ell,j_0}^{(k)}f_{1,k}
 +\sum_{k=1}^{d_2}(x_{2,k}+\mathrm{i}y_{2,k})Y_{\ell,j_0}^{(k)}f_{2,k}\right],\\
 \operatorname{Im}V&=\sqrt2\operatorname{Re}\left[
 \sum_{k=1}^{d_1}(x_{1,k}+\mathrm{i}y_{1,k})Y_{\ell,j_0}^{(k)}f_{1,k}
 -\sum_{k=1}^{d_2}(c_{2,k}+\mathrm{i}d_{2,k})Y_{\ell,j_0}^{(k)}f_{2,k}\right],
\end{aligned}
\end{equation}
where the real constants form coefficient vectors $c,x,d,y$, indexed by $(i,k)$.
The conditions $\langle V,e_1\rangle=\langle V,e_2\rangle=0$ become
\[
 \sqrt{\alpha_1}c_{1,1}-\sqrt{\alpha_2}c_{2,1}=0,\qquad
 \sqrt{\alpha_1}d_{1,1}-\sqrt{\alpha_2}d_{2,1}=0.
\]
Define $v$ by
\[
 v_{1,1}=\sqrt{\alpha_1},\qquad v_{2,1}=-\sqrt{\alpha_2},\qquad
 v_{i,k}=0\quad(k>1).
\]
Then these constraints read $c,d\perp v$.
The vectors $x,y$ are unrestricted.
Each vector has $d_1+d_2$ entries.

We set $T=T_1\oplus T_2$, where $T_i$ is the $d_i\times d_i$
tridiagonal matrix with entries
\[
 (T_i)_{kk}=L-2k^2,\qquad
 (T_i)_{k,k+1}=(T_i)_{k+1,k}
 =-\sqrt{(\alpha_i-k(k+1))(L-k(k+1))}.
\]
Since $\varpi=0$, the frame and scalar ladder formulas give
\[
 \sum_{a=1}^2 e_ae_a^t=\frac K2(\alpha_1-\alpha_2),\qquad
 \int_{S^2}|\dd V|^2\,\dd A
 =K\sum_{u\in\{c,x,d,y\}}
       \left(u^tTu+\sum_{i=1}^2\sum_{k=1}^{d_i}\alpha_i u_{i,k}^2\right).
\]
The real coordinates also give
\[
 \int_{S^2}\operatorname{Re}(VV^t)\,\dd A
 =\sum_{i=1}^2(-1)^{i-1}\sum_{k=1}^{d_i}
       \bigl(c_{i,k}^2+d_{i,k}^2-x_{i,k}^2-y_{i,k}^2\bigr).
\]
Using $\alpha_1+\alpha_2=4/K$, we obtain the shifted matrix for the
$x$- and $y$-blocks:
\[
 D=\operatorname{diag}\bigl(T_1+(\alpha_1-\alpha_2)I,
                           T_2-(\alpha_1-\alpha_2)I\bigr).
\]
\Needspace{7\baselineskip}
The remaining integrals are
\[
 \int_{S^2}\sum_{a=1}^2\langle V,\mathrm i e_a\rangle^2\,\dd A
 =\frac K4\bigl\{(v^tx)^2+(v^ty)^2\bigr\},
\]
\[
 \int_{S^2}\sum_{a=1}^2\langle V,\mathrm{II}_{1a}\rangle^2\,\dd A
 =\frac{K^2}{16}\sum_{u\in\{c,d\}}
 \left(\sqrt{\alpha_1(\alpha_1-2)}\,u_{1,2}
       -\sqrt{\alpha_2(\alpha_2-2)}\,u_{2,2}\right)^2.
\]
Define $A_\ell$ on $v^\perp$ and $B_\ell$ on the full coefficient space by
\[
 A_\ell(c)=c^tTc-\frac K4
 \left(\sqrt{\alpha_1(\alpha_1-2)}\,c_{1,2}
       -\sqrt{\alpha_2(\alpha_2-2)}\,c_{2,2}\right)^2
 \quad(c\perp v),
\qquad 
 B_\ell(x)=x^tDx-(v^tx)^2.
\]
By orthonormality and \eqref{common:index-form}, we obtain
\begin{equation}\label{matrix:zero-form}
 \int_{S^2}|V|^2\,\dd A=c^tc+x^tx+d^td+y^ty,\qquad
 Q_f(V,V)=K\{A_\ell(c)+B_\ell(x)+A_\ell(d)+B_\ell(y)\}.
\end{equation}

For the following matrix expressions, choose orthonormal bases of the
coefficient spaces, using the same basis for repeated copies.
We retain the same symbols for the coefficient vectors and for the
matrices of the quadratic forms.
Then \eqref{matrix:nonzero-form} and \eqref{matrix:zero-form} become, respectively,
\[
 Q_f(V,V)=K
 \begin{pmatrix}x^t&y^t\end{pmatrix}
 \begin{pmatrix}M_\ell&0\\0&M_\ell\end{pmatrix}
 \begin{pmatrix}x\\y\end{pmatrix},
\quad 
 Q_f(V,V)=K
 \begin{pmatrix}c^t&x^t&d^t&y^t\end{pmatrix}
 \begin{pmatrix}
 A_\ell&0&0&0\\
 0&B_\ell&0&0\\
 0&0&A_\ell&0\\
 0&0&0&B_\ell
 \end{pmatrix}
 \begin{pmatrix}c\\x\\d\\y\end{pmatrix}.
\]
The norm identities in \eqref{matrix:nonzero-form} and \eqref{matrix:zero-form}
show that the coordinate groups are mutually orthogonal in $L^2$.
Thus, by polarization, the identity
\[
 Q_f(V,W)=-\int_{S^2}\langle JV,W\rangle\,\dd A
\]
determines the block diagonal form of $J_\ell$ for $\ell\ge1$.
Its order follows from the normality constraints.

\noindent\begin{minipage}{\linewidth}
Collecting the preceding calculations, we obtain the following.

\begin{proposition}[Jacobi matrices]\label{common:explicit-Jacobi}
In real orthonormal coefficient coordinates, the matrices in
\eqref{common:Jacobi-decomposition} have the following forms.

\textup{(i)} For Types~I and II, $J_0=0_{2(c-1)}$, and, for $\ell\ge1$,
\begin{equation}\label{matrix:nonzero-J}
 J_\ell=-K\operatorname{diag}(M_\ell,M_\ell).
\end{equation}
Here $M_\ell$ represents this quadratic form in an orthonormal basis
of the real space defined by the constraints on $u$.
The order of $J_\ell$ is
$2\{c(p+r+1)-2-\mathbf1_{\{\ell\ge\lambda\}}\}$.
Here $\mathbf1_{\{\cdot\}}$ equals $1$ when the stated condition holds
and $0$ otherwise.

\textup{(ii)} For Types~III.I and III.II, there is no $\ell=0$ normal mode,
and, for $\ell\ge1$,
\begin{equation}\label{directIII:J}
 J_\ell=-K\operatorname{diag}(A_\ell,B_\ell,A_\ell,B_\ell).
\end{equation}
The blocks correspond to $(c,x,d,y)$, with $A_\ell$ represented in
orthonormal coordinates on $v^\perp$ and $B_\ell$ on the full coefficient
space. The order of $J_\ell$ is $4(d_1+d_2)-2$.
\end{proposition}
\end{minipage}\par\smallskip

The positive and zero eigenvalue counts of these matrices give the index
and nullity through \eqref{common:real-multiplicity-counts}.
The calculations in Appendix~\ref{app:congruence} yield the following formulas.

\begin{proposition}[Inertia of the Jacobi matrices]\label{index:inertia}
The positive and zero eigenvalue counts of $J_\ell$ are as follows.

\medskip\noindent\textup{(i) Types I and II.}
For $1\le\ell\le n_0$,
\[
\begin{aligned}
 \nu_+(J_\ell)&=2c\min\{\ell,(n_0-\lambda)/2,n_0-\ell\}
 -2\mathbf1_{\{\ell\ge\lambda,\;n_0+\ell\text{ odd}\}},\\
 \nu_0(J_\ell)&=2c-2\mathbf1_{\{\ell\ge\lambda\}}
 +4\mathbf1_{\{\ell\ge\lambda,\;n_0+\ell\text{ odd}\}}-2\mathbf1_{\{\ell=1\}}.
\end{aligned}
\]
For $\ell>n_0$, $\nu_+(J_\ell)=\nu_0(J_\ell)=0$.

\medskip\noindent\textup{(ii) Types III.I and III.II.}
For $\ell\ge1$,
\begin{equation}\label{zero:direct-counts}
\begin{aligned}
 \nu_+(J_\ell)&=2\sum_{i=1}^2\max\left\{0,d_i-\left\lfloor\frac{\ell+1}{2}\right\rfloor\right\}
       +2\max\{0,\min(\ell,m_2/2,n/2-\ell)\},\\
 \nu_0(J_\ell)&=2\mathbf1_{\{\ell\text{ odd}\}}\sum_{i=1}^2\mathbf1_{\{\ell<m_i\}}
         +2\mathbf1_{\{\ell\le n/2\}}-2\mathbf1_{\{\ell=1\}}.
\end{aligned}
\end{equation}
\end{proposition}

\Needspace{12\baselineskip}
\begin{proof}[Proof of Theorem~\ref{main:theorem}]
By \eqref{common:real-multiplicity-counts}, the index and nullity are obtained
by summing the counts in Proposition~\ref{index:inertia} with multiplicity $2\ell+1$.
For Types~I and II, the $\ell=0$ contribution to the nullity is
$2(c-1)$. Reordering the finite sums gives
\[
 \sum_{\ell=0}^{n_0}(2\ell+1)\min\{\ell,(n_0-\lambda)/2,n_0-\ell\}
 =\sum_{j=1}^{(n_0-\lambda)/2}\sum_{\ell=j}^{n_0-j}(2\ell+1)
 =\frac{(n_0^2-\lambda^2)(n_0+1)}4,
\]
\[
 \sum_{\substack{\lambda\le\ell\le n_0\\n_0+\ell\text{ odd}}}(2\ell+1)
 =\sum_{j=0}^{(n_0-\lambda)/2-1}\bigl[2(\lambda+1+2j)+1\bigr]
 =\frac{(n_0-\lambda)(n_0+\lambda+1)}2.
\]
Sums over empty index ranges are understood to be zero.
These identities yield
\begin{equation}\label{main:nonzero}
 \Ind f=\frac{n_0-\lambda}{2}\bigl[n(n_0+\lambda)-2\bigr],\qquad
 \Nul f=2n(n_0+1)+2(n_0-\lambda)-6.
\end{equation}

For Types~III.I and III.II, the same procedure gives, for $i=1,2$,
\[
 \sum_{\ell\ge1}(2\ell+1)
 \max\left\{0,\min\{m_i/2,\ell\}-\left\lfloor\frac{\ell+1}{2}\right\rfloor\right\}
 =\sum_{j=1}^{m_i/2}\sum_{\ell=j}^{2j-2}(2\ell+1)
 =\frac{m_i(m_i-2)(m_i+1)}8,
\]
\[
 \sum_{\ell\ge1}(2\ell+1)
 \max\{0,\min(\ell,m_2/2,n/2-\ell)\}
 =\sum_{j=1}^{m_2/2}\sum_{\ell=j}^{n/2-j}(2\ell+1)
 =\frac{m_1m_2(n+2)}8.
\]
Consequently,
\[
 \Ind f=\sum_{\ell\ge1}(2\ell+1)\nu_+(J_\ell)
 =\frac14\left\{\sum_{i=1}^2m_i(m_i-2)(m_i+1)+m_1m_2(n+2)\right\}.
\]
Summing the zero-eigenvalue counts in \eqref{zero:direct-counts} gives
\[
 \Nul f=\sum_{\ell\ge1}(2\ell+1)\nu_0(J_\ell)
 =\sum_{i=1}^2m_i(m_i+1)+\frac{n(n+4)}2-6.
\]
Since $n=m_1+m_2$, the curvature formula in
Theorem~\ref{hom:classification} reads
\[
 \frac{16}{K}=m_1^2+m_2^2+2n.
\]
Substitution gives
\begin{equation}\label{main:zero}
 \Ind f=(n-2)\left(\frac4K-\frac n4\right),\qquad
 \Nul f=\frac{16}{K}+\frac{n(n+2)}2-6.
\end{equation}
Equations~\eqref{main:nonzero} and \eqref{main:zero} give all four
rows of Table~\ref{models:formulas} and imply
\eqref{eq:main-unified}.
\end{proof}

For example, in Type~III.I with $m=2$, we have $K=1$ and, for $\ell\ge1$,
\[
 A_\ell=L-2,\qquad
 B_\ell=\begin{pmatrix}L-4&2\\2&L-4\end{pmatrix},
 \qquad L=\ell(\ell+1).
\]
Thus $J_\ell$ in \eqref{directIII:J} is orthogonally similar to
$\operatorname{diag}((2-L)I_4,(6-L)I_2)$.
The modes $\ell=1,2$ give, respectively,
$(\nu_+(J_\ell),\nu_0(J_\ell))=(2,4),(0,2)$;
all higher modes are negative definite. Hence
\[
 \Ind f=3\cdot2=6,\qquad \Nul f=3\cdot4+5\cdot2=22.
\]

\section{Integrability of normal Jacobi fields}\label{sec:integrability}

Theorem~\ref{main:integrability} is proved by deforming the model frames,
with Types~I and II treated in Section~\ref{sec:integrability-nonzero}
and Types~III.I and III.II in Section~\ref{sec:integrability-zero}.
Both constructions use the following dimension estimate.

\begin{lemma}\label{int:harmonic-dimension}
Let $f:S^2\to Q_n$ be a conformal minimal immersion, and let
$\mathcal V$ be a real vector space of fields along $f$ satisfying
the linearized harmonic map equation. The normal projection
$\mathcal V\to\ker J_f$, $V\mapsto V^\perp$, is well defined and has
kernel of dimension at most six. Consequently,
\begin{equation}\label{int:common-dimension-bound}
 \dim_\R\mathcal V\le\Nul f+6.
\end{equation}
If equality holds, the normal projection is surjective.
\end{lemma}
\begin{proof}
For $V\in\mathcal V$, let $f_t$ be a smooth variation satisfying
\[
 f_0=f,\qquad \left.\partial_t f_t\right|_{t=0}=V,
\]
and let $g_t$ be its induced metric, with $g=g_0$ and
$\dot g_0=\left.\partial_tg_t\right|_{t=0}$.
Since $V$ satisfies the linearized harmonic map equation, the linearized
Hopf differential is a holomorphic quadratic differential on $S^2$, and hence vanishes
\cite{WoodJacobi}. First-order conformality and the normal part of the
linearized harmonic map equation therefore give
\[
 \dot g_0=\frac12\bigl(\operatorname{tr}_g\dot g_0\bigr)g,
 \qquad J_f(V^\perp)=0.
\]
For $V^\perp=0$, the tangent vector field $\xi$ defined by $V=\dd f(\xi)$ satisfies
\[
 \mathcal L_\xi g=\dot g_0.
\]
Thus $\xi$ is conformal. The space of conformal vector fields on $S^2$
has real dimension six \cite{LemaireWood}, and hence
\[
 \dim_\R\mathcal V
 =\dim_\R\{V\in\mathcal V:V^\perp=0\}
 +\dim_\R\{V^\perp:V\in\mathcal V\}
 \le6+\Nul f.
\]
If equality holds, then
\[
 \dim_\R\{V^\perp:V\in\mathcal V\}=\Nul f,\qquad
 \{V^\perp:V\in\mathcal V\}=\ker J_f.\qedhere
\]
\end{proof}

\subsection{Types I and II}\label{sec:integrability-nonzero}
Let $f=[Z_\lambda]$ be a model of Type~I or II, with the frame $Z_\nu$
and fixed positive parameters $(n_0,\lambda)$ as in Section~\ref{sec:setup}.
Fix an orthonormal basis of $\C^{n+2}$ extending
$Z_{n_0}(1),Z_{n_0-2}(1),\ldots,Z_{-n_0}(1)$.
Let $\mathsf C$ be a constant matrix near $I_{n+2}$ satisfying
\[
 x\mathsf C=x\qquad\text{for every added basis vector }x.
\]
In the original ambient coordinates, set
\[
 c_\nu(g)=Z_\nu(1)\rho(g)\mathsf C=Z_\nu(g)\mathsf C.
\]
In the chosen basis, the coordinates of the initial vectors $c_\nu(1)$
form the free rows of $\mathsf C$.
Since $\mathsf C$ is constant, the vectors $c_\nu$ satisfy
\eqref{int:nonzero-ladder}.
Let $F=\operatorname{span}_\C\{c_\nu:\nu\ge\lambda\}$ and impose
\begin{equation}\label{int:nonzero-frame-equations}
 c_\nu(g)c_\nu(g)^t=0\qquad(\nu\ge\lambda,\ g\in\SU(2)).
\end{equation}
At $\mathsf C=I_{n+2}$, $c_\nu=Z_\nu$, so
\eqref{int:nonzero-frame-equations} holds.
Differentiating \eqref{int:nonzero-frame-equations} using
\eqref{int:nonzero-ladder} gives $c_\nu c_\mu^t=0$ for
$\nu+\mu\ge2\lambda-2$.
Conjugation gives $\bar F=\operatorname{span}_\C\{\bar c_\mu:\mu\ge\lambda\}$.
Let $S$ be an independent constant complex matrix satisfying
$S^t=-S$ and $xS=0$ for each added basis vector $x$.
Its independent entries are $Z_\nu(1)SZ_\mu(1)^t$,
$\nu>\mu$ in $\Delta_{n_0}$, chosen near zero.
For $\nu>\lambda$, define $u_\nu$ by
\begin{equation}\label{int:alternating-frame-correction}
 u_\nu c_\mu^t=Z_\nu(g)S Z_\mu(g)^t,\qquad
 u_\nu-c_\nu\in\bar F
 \quad(\mu\ge\lambda).
\end{equation}
Since $\{\bar c_\mu:\mu\ge\lambda\}$ is a basis of $\bar F$,
these pairings determine $u_\nu-c_\nu$ uniquely; in particular,
$u_\nu=c_\nu$ at $S=0$.
Since $S^t=-S$, \eqref{int:alternating-frame-correction} gives
\[
 u_\nu u_\mu^t
 =u_\nu c_\mu^t+u_\mu c_\nu^t=0
 \qquad(\nu,\mu>\lambda).
\]
Set $U=\operatorname{span}_\C\{u_\nu:\nu>\lambda\}$ and define $P$
by the Hermitian orthogonal decomposition
\begin{equation}\label{int:nonzero-frame-family}
 F\oplus\bar F=P\oplus U\oplus\bar U.
\end{equation}
Linear independence persists near the model, so
\[
 \dim_\C F-1=\dim_\C U=\frac{n_0-\lambda}{2},\qquad
 \dim_\C P=2(\dim_\C F-\dim_\C U)=2,\qquad \bar P=P.
\]
At the model, $P(g)=\C Z_\lambda(g)\oplus\C\bar Z_\lambda(g)$.
Let $Z(g)\in P(g)$ be a local unit vector such that $Z(g)Z(g)^t=0$
and $\C Z(g)$ is close to $\C Z_\lambda(g)$.
Then
\[
 P(g)=\C Z(g)\oplus\C\bar Z(g),\qquad
 \C^{n+2}=\C Z(g)\oplus\C\bar Z(g)\oplus P(g)^\perp.
\]
The frame transformation laws imply that $\C Z(g)$ depends only on $[g]$.
Thus the following map is well defined:
\[
 f_{\mathsf C,S}:S^2\to Q_n,\qquad
 f_{\mathsf C,S}([g])=[Z(g)],\qquad f_{I_{n+2},0}=f.
\]

\begin{lemma}\label{int:nonzero-frame-geometry}
Suppose that $\mathsf C$ satisfies \eqref{int:nonzero-frame-equations}.
For $(\mathsf C,S)$ sufficiently close to $(I_{n+2},0)$,
the map $f_{\mathsf C,S}$ is a conformal minimal immersion.
At the model, solutions of the linearized equations with
$\dot f_{\mathsf C,S}=0$ satisfy $\dot F=\dot U=0$.
Dots denote first variations at $(\mathsf C,S)=(I_{n+2},0)$.
\end{lemma}
\begin{proof}
Write $\dd u_\nu=(\dd u_\nu)_\omega\,\omega
+(\dd u_\nu)_\varphi\,\varphi
+(\dd u_\nu)_{\bar\varphi}\,\bar\varphi$,
where subscripts denote coframe coefficients.
Equation \eqref{int:nonzero-ladder} gives $(\dd c_\nu)_{\bar\varphi}\in F$
for $\nu\ge\lambda$.
Differentiation of \eqref{int:alternating-frame-correction}, with
$\mathsf C,S$ constant, gives
\[
 (\dd u_\nu)_\varphi\in F\oplus\bar F,\qquad
 (\dd u_\nu)_\varphi u_\mu^t=0,\qquad
 (\dd u_\nu)_{\bar\varphi}-a_{\nu+2,n_0}u_{\nu+2}
 \in(F\oplus\bar F)^\perp
 \quad(\nu,\mu>\lambda).
\]
Hence $(\dd u_\nu)_\varphi\in P\oplus U$ and
$(\dd u_\nu)_{\bar\varphi}\in U\oplus(F\oplus\bar F)^\perp$.
The line $P\cap(F+\bar U)$ is close to $\C Z_\lambda$.
Let $Z$ be a local unit section of this line and let $w\in F$ satisfy $Z-w\in\bar U$.
Then $w\perp U$, and
\[
 ZZ^t=2w(Z-w)^t=0.
\]
Since $(\dd w)_{\bar\varphi}\in F$ and
$(\dd\bar u_\nu)_{\bar\varphi}\in P\oplus\bar U$,
it follows that $(\dd Z)_{\bar\varphi}\in F\oplus\bar F$.
The decomposition of $\dd Z$ and differentiation of $Z\perp U,\bar U$ give
\[
 \dd Z=\varpi Z+X\varphi+Y\bar\varphi,\qquad
 X,Y\in P^\perp,\qquad Y\in U,\qquad X\perp U.
\]
By \eqref{common:tangent-covariant-definition},
$(\dd u_\nu)_\varphi\in P\oplus U$ and $\dd^2Z=0$ give
\begin{equation}\label{int:mixed-projection}
 (DX)_{\bar\varphi}\in U^\perp,\qquad
 (DY)_\varphi\in P\oplus U,\qquad
 (DX)_{\bar\varphi}-(DY)_\varphi\in P.
\end{equation}
Thus $(DX)_{\bar\varphi},(DY)_\varphi\in P$, so the map is harmonic.
Its Hopf differential is a holomorphic quadratic differential on $S^2$,
and hence vanishes. Thus the map is conformal and remains an immersion
near the model.

Suppose $\dot f_{\mathsf C,S}=0$, equivalently $\dot P=0$,
and $\lambda<n_0$.
The vectors $\dot c_\nu$ also satisfy \eqref{int:nonzero-ladder}.
Linearizing $u_\nu\perp P$ gives
$\langle\dot c_\nu,Z_\lambda\rangle_\C=0$ for $\nu>\lambda$.
If $\langle\dot c_\nu,Z_\mu\rangle_\C\equiv0$, differentiation gives
\[
a_{\nu,n_0}\langle\dot c_{\nu-2},Z_\mu\rangle_\C
 =a_{\mu+2,n_0}\langle\dot c_\nu,Z_{\mu+2}\rangle_\C,\qquad
a_{\mu,n_0}\langle\dot c_\nu,Z_{\mu-2}\rangle_\C
 =a_{\nu+2,n_0}\langle\dot c_{\nu+2},Z_\mu\rangle_\C.
\]
Terms beyond the weight chain are omitted.
Induction in the first relation from $\mu=\lambda$, followed by the second
relation at $\nu=n_0$, gives
\[
 \langle\dot c_\nu,Z_\mu\rangle_\C=0\quad(\nu>\mu\ge\lambda),\qquad
 \langle\dot c_{n_0},Z_\mu\rangle_\C=0\quad(\mu<n_0).
\]
Set $\gamma=\langle\dot c_{n_0},Z_{n_0}\rangle_\C$.
Applying the first relation down the weight chain gives
$\langle\dot c_\nu,Z_\nu\rangle_\C=\gamma$ for all $\nu\in\Delta_{n_0}$.
The adjacent terms in \eqref{int:nonzero-ladder} cancel in the sum
\[
 (n_0+1)\dd\gamma
 =\dd\sum_{\nu\in\Delta_{n_0}}\langle\dot c_\nu,Z_\nu\rangle_\C=0.
\]
Thus $\gamma$ is constant. Comparing $\varphi$ coefficients successively gives
\[
 \langle\dot c_\nu,Z_\mu\rangle_\C=\gamma\delta_{\nu\mu}
 \qquad(\nu,\mu\in\Delta_{n_0}).
\]
After the replacement $\dot c_\nu\mapsto\dot c_\nu-\gamma Z_\nu$,
$\langle\dot c_\nu,Z_\mu\rangle_\C=0$ for all $\nu,\mu\in\Delta_{n_0}$.
Linearizing $P\subset F\oplus\bar F$ and $c_\lambda c_\lambda^t=0$ gives
$\dot c_\lambda\in\operatorname{span}_\C\{\bar Z_\nu:\nu>\lambda\}$.
The ladder relations and the preceding inclusion give
\[
 \dot c_\nu\in\operatorname{span}_\C\{\bar Z_\mu:\mu\in\Delta_{n_0}\}
 \quad(\nu\in\Delta_{n_0}),\qquad
 \dot c_\lambda Z_\mu^t=0\quad(\mu\le\lambda).
\]
If $\dot c_\nu Z_\mu^t\equiv0$, differentiation gives
\[
a_{\nu,n_0}\dot c_{\nu-2}Z_\mu^t
 +a_{\mu,n_0}\dot c_\nu Z_{\mu-2}^t=0,\qquad
a_{\nu+2,n_0}\dot c_{\nu+2}Z_\mu^t
 +a_{\mu+2,n_0}\dot c_\nu Z_{\mu+2}^t=0.
\]
Successive applications of the second recurrence give
\[
 \dot c_{\lambda+2j}Z_\mu^t=0
 \quad\left(0\le j\le\frac{n_0-\lambda}{2},\quad
 \mu\in\Delta_{n_0},\quad\mu\le\lambda-2j\right).
\]
Using the first recurrence at $\mu=-n_0$, where $a_{-n_0,n_0}=0$,
and then the second recurrence yields
\[
 \dot c_\lambda Z_\mu^t=0\ (\mu\le\lambda)
 \;\Longrightarrow\;
 \dot c_\nu Z_{-n_0}^t=0\ (\nu\ge\lambda)
 \;\Longrightarrow\;
 \dot c_\nu Z_{-n_0}^t=0\ (\nu\in\Delta_{n_0})
 \;\Longrightarrow\;
 \dot c_\nu Z_\mu^t=0\ (\nu,\mu\in\Delta_{n_0}).
\]
Since $\dot c_\nu$ lies in the span of the $\bar Z_\mu$, it follows that
$\dot c_\nu=0$ for every $\nu\in\Delta_{n_0}$.
Since $\dot P=0$ and $\dot c_\nu=0$, differentiating $u_\nu\perp P$
and \eqref{int:alternating-frame-correction} gives
$Z_\nu(g)\dot S Z_\lambda(g)^t=0$ for $\nu>\lambda$.
Applying the raising relation successively in the second index yields
\[
 Z_\nu(g)\dot S Z_\mu(g)^t=0
 \qquad(\nu>\lambda,\ \mu\ge\lambda,\ g\in\SU(2)).
\]
If $\lambda=n_0$, linearizing $c_{n_0}c_{n_0}^t=0$ with $\dot P=0$
gives $\dot c_{n_0}\in\C Z_{n_0}$.
The same diagonal-pairing argument shows that $\gamma$ is constant.
Comparing $\varphi$ coefficients in \eqref{int:nonzero-ladder} then gives
$\dot c_\nu=\gamma Z_\nu$ for all $\nu$; in this case $U=0$.

Restoring the scalar term establishes, in both cases, the forward implication in
\begin{equation}\label{int:nonzero-zero-field}
 \dot P=0\quad\Longleftrightarrow\quad
 \begin{cases}
  \dot c_\nu=\gamma Z_\nu,&\nu\in\Delta_{n_0},\\
  Z_\nu(g)\dot S Z_\mu(g)^t=0,&\nu>\lambda,\ \mu\ge\lambda,\ \forall g,
 \end{cases}
\end{equation}
where $\gamma\in\C$ is constant. Conversely,
\eqref{int:alternating-frame-correction} gives $\dot u_\nu=\gamma Z_\nu$
for $\nu>\lambda$, hence $\dot F=\dot U=\dot P=0$.
\end{proof}

\begin{proof}[Proof of Theorem~\ref{main:integrability} (for Types~I and II)]
The simultaneous rescaling preserves $P$:
\[
 (c_\nu(1),S)\longmapsto(c\,c_\nu(1),c^2S)\qquad(c\ne0).
\]
The parameters can therefore be normalized by
\[
 \langle c_{n_0}(1),Z_{n_0}(1)\rangle_\C=1.
\]
Within the space of matrices satisfying $S^t=-S$ and $xS=0$ for each
added basis vector $x$, let $S$ range over a fixed linear complement of
the subspace defined by
\[
 Z_\nu(g)S Z_\mu(g)^t=0
 \qquad(\nu>\lambda,\ \mu\ge\lambda,\ g\in\SU(2)).
\]
For normalized linearized solutions, \eqref{int:nonzero-zero-field} gives
\[
 \dot P=0\quad\Longrightarrow\quad\dot{\mathsf C}=\dot S=0.
\]

The coefficients $c_\nu(1)c_\mu(1)^t$ are symmetric in $\nu,\mu$,
whereas $Z_\nu(1)S Z_\mu(1)^t$ are alternating; their independent
entries are indexed by $\nu\le\mu$ and $\nu<\mu$, respectively.
For each integer $0\le r\le n_0$, the numbers of coefficients of total weight $2r$ are
\[
 \begin{aligned}
 \#\{(\nu,\mu)\in\Delta_{n_0}^{\,2}:\nu\le\mu,\ \nu+\mu=2r\}
 &=\left\lfloor\frac{n_0-r}{2}\right\rfloor+1,\\
 \#\{(\nu,\mu)\in\Delta_{n_0}^{\,2}:\nu<\mu,\ \nu+\mu=2r\}
 &=\left\lfloor\frac{n_0-r+1}{2}\right\rfloor.
 \end{aligned}
\]
Differences of successive weight multiplicities give highest weights
$2n_0,2n_0-4,\ldots$ in the symmetric case and
$2n_0-2,2n_0-6,\ldots$ in the alternating case, each with multiplicity one.
For a summand of highest weight $m$, let $h_{\nu\mu}$ be the
coefficients of a nonzero highest-weight vector in the product basis
$v_{\nu,n_0}\otimes v_{\mu,n_0}$, where $\nu+\mu=m$.
The raising relation on the two factors gives
\[
 a_{\nu+2,n_0}h_{\nu\mu}
 +a_{\mu,n_0}h_{\nu+2,\mu-2}=0.
\]
Whenever both index pairs lie in $\Delta_{n_0}^2$, the two ladder
coefficients are nonzero. This recurrence gives $h_{\nu\mu}\ne0$
for every allowed pair, with independent pairs indexed by
$\nu\le\mu$ in the symmetric case and $\nu<\mu$ in the alternating case.
In particular, the symmetric coefficient with $\nu=\mu=m/2$ is nonzero.
Hence \eqref{int:nonzero-frame-equations} eliminates precisely the
symmetric summands with highest weights $2\lambda,2\lambda+4,\ldots,2n_0$.
Since $\dim_\C V_m=m+1$, the number of complex coefficient equations is
\[
 \sum_{j=0}^{(n_0-\lambda)/2}(2\lambda+4j+1)
 =\frac{(n_0-\lambda+2)(n_0+\lambda+1)}2.
\]
These equations are linear in the pairings $c_\nu(1)c_\mu(1)^t$.
Their linearizations are therefore the same linear relations applied to
\[
 \left.\partial_t\right|_{t=0}\bigl(c_\nu(1)c_\mu(1)^t\bigr)
 =\dot c_\nu(1)Z_\mu(1)^t+Z_\nu(1)\dot c_\mu(1)^t.
\]
In the alternating case, skew-symmetry also gives a nonzero coefficient at
$(\nu,\mu)=(m/2+1,m/2-1)$.
By \eqref{int:nonzero-zero-field}, the number of independent complex
$S$-parameters after normalization is the total dimension of the
alternating summands with highest weights
$2\lambda+2,2\lambda+6,\ldots,2n_0-2$:
\[
 \sum_{j=0}^{(n_0-\lambda)/2-1}(2\lambda+4j+3)
 =\frac{(n_0-\lambda)(n_0+\lambda+1)}2,
\]
with an empty sum when $\lambda=n_0$.
The preceding counts give the difference between the numbers of
normalized real parameters and real equations:
\begin{equation}\label{int:nonzero-parameter-count}
 2(n+2)(n_0+1)-2-2(n_0+\lambda+1)
 =2n(n_0+1)+2n_0-2\lambda=\Nul f+6.
\end{equation}
The last equality follows from \eqref{main:nonzero}.
Let $\mathcal V$ be the real vector space of fields induced by normalized
linearized solutions. Linearizing the identities in the proof of
Lemma~\ref{int:nonzero-frame-geometry} shows that these fields satisfy the
linearized harmonic map equation. Equation~\eqref{int:nonzero-zero-field}
shows that the map from normalized linearized solutions to their induced
fields is injective.
Equation~\eqref{int:nonzero-parameter-count} and
Lemma~\ref{int:harmonic-dimension} give
\[
 \Nul f+6\le\dim_\R\mathcal V\le\Nul f+6,\qquad
 \{V^\perp:V\in\mathcal V\}=\ker J_f.
\]
Thus the linearized coefficient system has full row rank:
\[
 \operatorname{rank}_{\R}(\text{linearized system})
 =(n_0-\lambda+2)(n_0+\lambda+1).
\]
The implicit function theorem and Lemma~\ref{int:nonzero-frame-geometry}
give a local family of conformal minimal immersions whose variation fields
at the model form $\mathcal V$. Choose a real subspace $E\subset\mathcal V$ with
\[
 \mathcal V=\{V\in\mathcal V:V^\perp=0\}\oplus E,\qquad
 \{V^\perp:V\in E\}=\ker J_f.
\]
Choose a local parametrization by a neighbourhood of $0$ in $\mathcal V$,
with differential equal to the identity at $0$.
Restricting it to $E$ and using the linear isomorphism
$E\to\ker J_f$, $V\mapsto V^\perp$, gives a subfamily
parametrized by $v$ near $0$ in $\ker J_f$ such that
\[
 \Phi(0,\cdot)=f,\qquad
 \left(\left.\partial_t\Phi(tv,\cdot)\right|_{t=0}\right)^\perp=v
 \quad(v\in\ker J_f).
\]
\end{proof}

\subsection{Types III.I and III.II}\label{sec:integrability-zero}
For the frames in \eqref{zero:weight-frames}, with $m_1=m_2$ in
Type~III.I, let
\[
 F_0=\operatorname{span}_{\C}\{f_{i,a}:i=1,2,\ a>0\}.
\]
Order these positive-weight frame vectors by nonincreasing weight,
with a fixed order among equal weights. Along a smooth curve in the
parameter space through the model, denote the deformed frame vectors by
$w_j(t,g)$, where $g\in\SU(2)$ and $t=0$ corresponds to the model.
If the $j$th original vector is $f_{i,a}$, then
\[
 w_j(0,g)=f_{i,a}(g),\qquad \lambda_j=2a,\qquad
 \dot w_j(g)=\left.\partial_t w_j(t,g)\right|_{t=0}.
\]
Dots denote derivatives at $t=0$ with $g$ fixed.
The differential $\dd$ and the operators $A,B$ act only on $g$;
the forms $\omega,\varphi,\bar\varphi$ are independent of the deformation
parameters. Pullback along fixed local sections of $\SU(2)\to S^2$
gives local frames, denoted by the same symbols.
The arguments $t,g$ are suppressed below unless needed.
Assume that the vectors $w_1,\ldots,w_{n/2}$ are linearly independent and that their spans
$F=\operatorname{span}_{\C}\{w_1,\ldots,w_{n/2}\}$ agree on overlaps.
Write
$\dd w_i=(\dd w_i)_\varphi\varphi+(\dd w_i)_{\bar\varphi}\bar\varphi$.
The frame equations are
\begin{equation}\label{zero:frame-equations}
 w_iw_j^t=0,\qquad (\dd w_i)_{\bar\varphi}\in F,\qquad
 (\dd w_i)_\varphi w_j^t=0\qquad(1\le i,j\le n/2).
\end{equation}
For $P=(F+\bar F)^\perp$, these equations give the Hermitian orthogonal decomposition
\[
 \C^{n+2}=P\oplus F\oplus\bar F,\qquad
 \dim_\C P=2,\qquad \bar P=P.
\]
At the model, $P=\operatorname{span}_{\C}\{f_{1,0},f_{2,0}\}$.
For $F$ close to $F_0$, the null line in $P(x)$ close to $f(x)$ has
a local unit lift $Z(x)$ satisfying
\[
 ZZ^t=0,\qquad Z\bar Z^t=1,\qquad P=\C Z\oplus\C\bar Z.
\]
These lines agree on overlaps and define
$f_F:S^2\to Q_n$, $f_F(x)=[Z(x)]$, with $f_{F_0}=f$.

\begin{lemma}\label{int:zero-frame-geometry}
Suppose that the local frames $w_1,\ldots,w_{n/2}$ satisfy
\eqref{zero:frame-equations}.
For $F$ sufficiently close to $F_0$, the map $f_F$ is a conformal
minimal immersion. At the model, solutions of the linearized frame equations
induce fields satisfying the linearized harmonic map equation, and
$\dot f_F=0$ implies $\dot F=0$.
\end{lemma}
\begin{proof}
The frame equations and $Z\perp F,\bar F$ give
\[
 \dd Z=\varpi Z+X\varphi+Y\bar\varphi,
 \quad X\in\bar F,\quad Y\in F,\qquad
 (\dd w_i)_\varphi\in P\oplus F\quad(1\le i\le n/2).
\]
Equation~\eqref{common:tangent-covariant-definition} and $\dd^2Z=0$ yield
\[
 (DX)_{\bar\varphi}\in P\oplus\bar F,\qquad
 (DY)_\varphi\in P\oplus F,\qquad
 (DX)_{\bar\varphi}-(DY)_\varphi\in P.
\]
Hence $(DX)_{\bar\varphi},(DY)_\varphi\in P$, so $f_F$ is harmonic.
The identity $X\bar Y^t=0$ gives conformality, and the immersion condition
persists near the model. Thus $f_F$ is a conformal minimal immersion.
For solutions of the linearized frame equations, linearization of these
identities gives the linearized harmonic map equation for every induced field.

For a linearized solution with zero induced field,
\[
 \dot f_F=0\quad\Longleftrightarrow\quad\dot P=0
 \quad\Longrightarrow\quad\dot w_j f_{i,0}^{\,t}=0
 \quad(i=1,2,\ 1\le j\le n/2).
\]
Fix $i\in\{1,2\}$ and $0\le a<m_i/2$.
If $\dot w_j f_{i,a}^{\,t}=0$ for every $j$, linearization of
$(\dd w_j)_{\bar\varphi}\in F$ and comparison of the
$\bar\varphi$ coefficients in $\dd(\dot w_j f_{i,a}^{\,t})=0$ give
\[
 a_{2a+2,m_i}\dot w_j f_{i,a+1}^{\,t}=0\qquad(1\le j\le n/2).
\]
Since $a_{2a+2,m_i}\ne0$, induction gives
$\dot w_j f_{i,a}^{\,t}=0$ for $i=1,2$, $0\le a\le m_i/2$ and
every $j$. The model frame decomposition therefore gives
\[
 \dot w_j\in F_0\quad(1\le j\le n/2),\qquad \dot F=0.
\]
\end{proof}

\begin{proof}[Proof of Theorem~\ref{main:integrability} (for Types~III.I and III.II)]
The weights fixed above satisfy
\begin{equation}\label{zero:weight-order}
 \lambda_1\ge\cdots\ge\lambda_{n/2}>0.
\end{equation}
The deformation preserves the weight $\lambda_i$ of each $w_i$:
\[
 w_i(t,\exp(s\varepsilon_1)g)=e^{\mathrm i\lambda_i s}w_i(t,g).
\]
For each fixed even weight $d$, let $[h]_{V_m}$ denote the projection of $h$
onto the span of the entries of $v_{d,m}\rho_m(g)$, where
$m=|d|,|d|+2,\ldots$. All functions below are finite sums of
such components; operators and projections act on each ambient coordinate
of a vector. The operators $A$ and $B$ decrease and increase the weight
by $2$, respectively.

The original frame relations are
\[
 f_{\alpha,a}B=-a_{2a+2,m_\alpha}f_{\alpha,a+1}
 \quad(0<a<m_\alpha/2),\qquad f_{\alpha,m_\alpha/2}B=0
 \quad(\alpha=1,2).
\]
After the fixed ordering, their coefficients are the model values
$L_{ij}(0)$; these are zero between the two original weight chains.
For $\lambda_j>\lambda_i$, choose $L_{ij}$ with the required weight
and the smallest admissible $m$:
\[
 L_{ij}\in V_{\lambda_j-\lambda_i-2}
 \quad\text{of weight }\lambda_i-\lambda_j+2\quad(\lambda_j>\lambda_i),
 \qquad L_{ij}=0\quad(\lambda_j\le\lambda_i).
\]
The free complex parameters are the constant coefficients of the functions
$L_{ij}$ and of each ambient coordinate of $[w_i]_{V_{\lambda_i}}$.
They vary independently near their model values. The variable $t$
parametrizes an arbitrary smooth curve through these values.

For functions of weight $d\ge0$, \eqref{common:scalar-mode-ladders} gives
\[
 \ker B=V_d,\qquad
 \operatorname{im}B=\{v:v\text{ has weight }d+2\}.
\]
For each fixed parameter value, solve successively
\begin{equation}\label{zero:triangular-frame}
 w_i(t,\cdot)B=\sum_{j<i}L_{ij}(t,\cdot)w_j(t,\cdot)
 \qquad(1\le i\le n/2)
\end{equation}
with the prescribed components $[w_i]_{V_{\lambda_i}}$.
At step $i$, the vectors $w_1,\ldots,w_{i-1}$ are already known, and
each term on the right has weight
$(\lambda_i-\lambda_j+2)+\lambda_j=\lambda_i+2$.
More explicitly,
\[
 w_1=[w_1]_{V_{\lambda_1}},\qquad
 [w_i]_{V_m}B=\left[\sum_{j<i}L_{ij}w_j\right]_{V_m}
 \quad(m>\lambda_i).
\]
Here $B$ changes the weight from $\lambda_i$ to $\lambda_i+2$ without
changing $m$.
For $m>\lambda_i$, the ladder coefficient in
\eqref{common:scalar-mode-ladders} is nonzero, so the second equation
determines $[w_i]_{V_m}$ uniquely. The remaining component
$[w_i]_{V_{\lambda_i}}$ is precisely the prescribed free component.
At the model parameter values, uniqueness recovers the original frame.

Products of matrix coefficients decompose into finitely many
representation spaces. Consequently, all $w_i$ lie in a fixed finite sum
of these representation spaces. Since each step uses only addition,
multiplication, and division by fixed nonzero ladder coefficients, the
resulting coefficients are polynomial, hence smooth, functions of the
free parameters.
With $F(t,g)=\operatorname{span}_{\C}\{w_1(t,g),\ldots,w_{n/2}(t,g)\}$,
the fixed weights give
\[
 F(t,\exp(s\varepsilon_1)g)=F(t,g).
\]
Thus $F$ is defined on $S^2$. Near the model, the vectors remain linearly
independent and $F\cap\bar F=0$, so $P=(F+\bar F)^\perp$ defines $f_F$
as above.

At the model, the parameter variations with $\dot F=0$ satisfy
\[
 \dot w_i=\sum_jG_{ij}w_j,\qquad G_{ij}=\dot w_i\bar w_j^t.
\]
All undifferentiated $w_i$ and $L_{ij}$ in the following linearized
formulas are evaluated at the model. Each $G_{ij}$ is a finite sum of
matrix coefficients of weight $\lambda_i-\lambda_j$, and
\eqref{zero:triangular-frame} gives
\[
 G_{ij}B=\sum_{k<i}L_{ik}G_{kj}
          -\sum_{k>j}G_{ik}L_{kj}+\dot L_{ij}.
\]
The equations for $G_{ij}$ are solved in increasing order of $i$ and, for each $i$,
decreasing order of $j$.
If $\lambda_i\ge\lambda_j$, then $\dot L_{ij}=0$ and the kernel
$V_{\lambda_i-\lambda_j}$ contributes $\lambda_i-\lambda_j+1$ free
complex coefficients. If $\lambda_i<\lambda_j$, then $G_{ij}B$ has no
$V_{\lambda_j-\lambda_i-2}$ component. Consequently,
\[
 \dot L_{ij}=-\left[
 \sum_{k<i}L_{ik}G_{kj}-\sum_{k>j}G_{ik}L_{kj}
 \right]_{V_{\lambda_j-\lambda_i-2}}.
\]
For every $m\ge\lambda_j-\lambda_i$, the ladder coefficient of $B$ in
\eqref{common:scalar-mode-ladders} is nonzero, so the equation determines
$[G_{ij}]_{V_m}$ uniquely.
Conversely, every choice of the free kernel coefficients determines all
$G_{ij}$ and $\dot L_{ij}$ by this recursion.
Together, $\dot L_{ij}$ and $[\dot w_i]_{V_{\lambda_i}}$, with
$\dot w_i=\sum_jG_{ij}w_j$, specify a parameter variation.
By uniqueness for the linearization of \eqref{zero:triangular-frame},
this variation induces precisely the vectors $\dot w_i$.

With the model as origin, restrict the parameters to a fixed complex
linear complement of these directions.
For each $(i,j)$, the number of $L_{ij}$ coefficients minus the number
of free $G_{ij}$ coefficients is $\lambda_j-\lambda_i-1$.
The remaining number of complex parameters is therefore
\[
 (n+2)\sum_i(\lambda_i+1)
 +\sum_{i,j=1}^{n/2}(\lambda_j-\lambda_i-1)
 =(n+2)\sum_i(\lambda_i+1)-\frac{n^2}{4}.
\]

For the deformed vectors, \eqref{zero:triangular-frame} gives
$(\dd w_i)_{\bar\varphi}\in F$ after pullback to local sections.
If $w_iw_j^t=0$ for all $i,j$, the $\omega$ term vanishes after pairing
with $w_j$, so $(\dd w_i)_\varphi w_j^t=(w_iA)w_j^t$.
Thus \eqref{zero:frame-equations} reduces to
$w_iw_j^t=(w_iA)w_j^t=0$ for all $i,j$.
Impose the finite coefficient equations
\begin{equation}\label{zero:finite-constraints}
 [w_iw_j^t]_{V_{\lambda_i+\lambda_j}}=0\quad(i\le j),\qquad
 [(w_iA)w_j^t]_{V_{\lambda_i+\lambda_j-2}}=0\quad(i<j).
\end{equation}
The product rule gives
\[
 (w_iw_j^t)B
 =\sum_{k<i}L_{ik}w_kw_j^t+\sum_{k<j}L_{jk}w_iw_k^t.
\]
Every pairing on the right has a smaller sum of indices.
Induction on $i+j$ therefore reduces $w_iw_j^t$ to its
$V_{\lambda_i+\lambda_j}$ component, which vanishes by
\eqref{zero:finite-constraints}.
Differentiation by $A$ then gives
$(w_iA)w_j^t=-(w_jA)w_i^t$.
Using $w_i(AB-BA)=\lambda_iw_i$ yields
\[
 \bigl((w_iA)w_j^t\bigr)B
 =\sum_{k<i}L_{ik}(w_kA)w_j^t
  +\sum_{k<j}L_{jk}(w_iA)w_k^t.
\]
The same induction at weight $\lambda_i+\lambda_j-2\ge0$ gives
$(w_iA)w_j^t=0$.
Hence \eqref{zero:finite-constraints} is equivalent to
\eqref{zero:frame-equations}.
At the model, $w_iw_j^t=(w_iA)w_j^t=0$.
For a solution of the linearization of \eqref{zero:finite-constraints},
linearizing the first recursion gives
$\partial_t(w_iw_j^t)|_{t=0}=0$.
The terms omitted from the second recursion are linear combinations
of the pairings $w_kw_l^t$; these pairings and their first variations
vanish at the model. The same recursion therefore proves that the
linearization of \eqref{zero:finite-constraints} is equivalent to the
linearized frame equations.
The number of complex coefficient equations is
\[
 \sum_{i\le j}(\lambda_i+\lambda_j+1)
 +\sum_{i<j}(\lambda_i+\lambda_j-1)
 =n\sum_i\lambda_i+\frac n2.
\]

Let $\mathcal V$ be the real vector space of fields induced by linearized
solutions. By Lemma~\ref{int:zero-frame-geometry} and the chosen complement,
these fields satisfy the linearized harmonic map equation, and the map
from linearized solutions to $\mathcal V$ is injective.
The positive weights in the two original chains are
$2,4,\ldots,m_\alpha$ for $\alpha=1,2$.
The preceding counts and the nullity computation in
Section~\ref{sec:index} give
\begin{equation}\label{zero:frame-dimension}
 4\sum_i\lambda_i+\frac{n(n+2)}2
 =\sum_{\alpha=1}^2m_\alpha(m_\alpha+2)+\frac{n(n+2)}2
 =\Nul f+6.
\end{equation}
Together with Lemma~\ref{int:harmonic-dimension}, this yields
\[
 \Nul f+6\le\dim_\R\mathcal V\le\Nul f+6.
\]
Hence the linearized coefficient system has full row rank, and equality
in Lemma~\ref{int:harmonic-dimension} gives
\[
 \operatorname{rank}_{\R}(\text{linearized system})=2n\sum_i\lambda_i+n,
 \qquad \{V^\perp:V\in\mathcal V\}=\ker J_f.
\]
The implicit function theorem and Lemma~\ref{int:zero-frame-geometry}
give a local family of conformal minimal immersions whose variation fields
at the model form $\mathcal V$. Choose a real subspace $E\subset\mathcal V$ with
\[
 \mathcal V=\{V\in\mathcal V:V^\perp=0\}\oplus E,\qquad
 \{V^\perp:V\in E\}=\ker J_f.
\]
Choose a local parametrization by a neighbourhood of $0$ in $\mathcal V$,
with differential equal to the identity at $0$.
Restricting it to $E$ and using the linear isomorphism
$E\to\ker J_f$, $V\mapsto V^\perp$, gives a subfamily
parametrized by $v$ near $0$ in $\ker J_f$ such that
\[
 \Phi(0,\cdot)=f,\qquad
 \left(\left.\partial_t\Phi(tv,\cdot)\right|_{t=0}\right)^\perp=v
 \quad(v\in\ker J_f).
\]
\end{proof}

\appendix
\titleformat{\section}{\normalfont\large\bfseries}{Appendix~\thesection}{1em}{}
\section{Inertia of the Jacobi matrices}\label{app:congruence}
The notation and matrices are those of Section~\ref{sec:jacobi-matrices}.

\subsection{Proof of Proposition~\ref{index:inertia}\textup{(i)}}\label{app:nonzero-congruence}

We determine the inertia of $M_\ell$ by restriction and completion of squares.
The full matrix of $q_\ell$ has eigenvalues
\begin{equation}\label{direct:scalar-spectrum}
 (\ell-j)(n_0+\ell-j+1),\qquad 0\le j\le p+r.
\end{equation}
To verify this, we set
$\omega_q=\binom{2\ell}{\ell+q}
\binom{n_0}{(n_0-\lambda)/2-q}$.
The change of coordinates $u_q=(-1)^q\sqrt{\omega_q}f(q)$ gives a
diagonally similar operator whose value on $f$ at $q$ is
\[
 \ell(n_0+\ell+1)f(q)
 +(\ell-q)\left(\frac{n_0-\lambda}{2}-q\right)
       \{f(q+1)-f(q)\}
 +(\ell+q)\left(\frac{n_0+\lambda}{2}+q\right)
       \{f(q-1)-f(q)\}.
\]
The endpoint coefficients vanish. Applied to $q^j$, this operator has
degree at most $j$, and the coefficient of $q^j$ is
\[
 \ell(n_0+\ell+1)-j(n_0+2\ell+1-j)
 =(\ell-j)(n_0+\ell-j+1).
\]
The evaluation vectors of
$1,q,\ldots,q^{p+r}$ form a basis, in which this operator is
triangular with the distinct diagonal entries
\eqref{direct:scalar-spectrum}.
Formula \eqref{direct:scalar-spectrum} also holds for $\lambda=0$
and even $n_0\ge0$.
For $1\le\ell\le n_0$, there are $\ell$ positive eigenvalues,
$p+r-\ell$ negative eigenvalues and one zero eigenvalue, where
$p+r-\ell=\min\{\ell,(n_0-\lambda)/2,n_0-\ell\}$.
For $\ell>n_0$, all eigenvalues are positive.

For $\ell\ge\lambda$, the constraint $VZ^t=0$ reduces to
$x_{-\lambda}=0$ in one coefficient block. In Type~II with $\tau>0$,
this uses the orthogonal change of coordinates
\[
 (x^{(1)},x^{(2)})\longmapsto
 \bigl(\cos(2\tau)x^{(1)}+\sin(2\tau)x^{(2)},
       -\sin(2\tau)x^{(1)}+\cos(2\tau)x^{(2)}\bigr),
\]
which preserves $q_\ell(x^{(1)})+q_\ell(x^{(2)})$.
Setting $x_{-\lambda}=0$ removes every term containing this coordinate,
so the matrix of the restricted form is obtained by deleting its
row and column.

For $\lambda\le\ell\le n_0$, we show that a nonzero kernel vector $w$
satisfies $w_{-\lambda}=0$ exactly when $n_0+\ell$ is odd.
With binomial coefficients outside their usual ranges taken as zero, we set
\[
 w_q=\frac1{\sqrt{\omega_q}}
 \sum_{a=0}^{\ell}(-1)^a\binom\ell a\binom\ell{a-q}
             \binom{n_0-\ell}{(n_0-\lambda)/2-a}.
\]
Applying the identity
\[
 \binom{m}{k}=\binom{m-1}{k-1}+\binom{m-1}{k}\qquad(m\ge1)
\]
to the binomial sum and shifting the summation index gives
$d_qw_q-b_{q-1}w_{q-1}-b_qw_{q+1}=0$, including the endpoints.
At $q=-p$, only the term $a=\ell-p$ survives, so $w_{-p}\ne0$.
Hence $w$ spans the kernel of the matrix of $q_\ell$.
Replacing $a$ by $\ell-\lambda-a$ in the sum gives
$w_{-\lambda}=(-1)^{\ell-\lambda}w_{-\lambda}$.
If $n_0+\ell$ is even, write $\ell=\lambda+2s$ and $n_0=\ell+2h$
with $s,h\in\mathbb Z_{\ge0}$, and set
\[
 t_{h,k}=(-1)^{s+k}\binom\ell{s+k}\binom\ell{s-k}\binom{2h}{h-k},
 \qquad
 G_{h,k}=\frac{(s+k)(\lambda+s+k)(h+1+k)}{2(h+1)}t_{h+1,k}.
\]
Then $\sqrt{\omega_{-\lambda}}w_{-\lambda}=\sum_{k=-s}^{s}t_{h,k}$,
and cancellation of factorials gives
\[
 (s+h+1)(\lambda+s+h+1)t_{h+1,k}
 -2(2h+1)(\ell+h+1)t_{h,k}=G_{h,k}-G_{h,k+1}.
\]
Since $G_{h,-s}=G_{h,s+1}=0$, summing gives
\[
 \sum_{k=-s}^{s}t_{h+1,k}
 =\frac{2(2h+1)(\ell+h+1)}{(s+h+1)(\lambda+s+h+1)}
   \sum_{k=-s}^{s}t_{h,k},\qquad
 \sum_{k=-s}^{s}t_{0,k}=(-1)^s\binom\ell s^2\ne0.
\]
The multiplier is positive, so the sum is nonzero for every $h\ge0$.

If $n_0+\ell$ is even, then $w_{-\lambda}\ne0$. The decomposition
$x=(x-x_{-\lambda}w/w_{-\lambda})+x_{-\lambda}w/w_{-\lambda}$
splits the coefficient space into the hyperplane $x_{-\lambda}=0$ and the kernel.

If $n_0+\ell$ is odd, then $-\lambda$ is an interior index.
Deleting this row and column gives two tridiagonal blocks.
By the endpoint recurrence, their kernels are spanned by the nonzero
left and right parts of $w$, respectively. The restricted form therefore has nullity two.
A codimension-one restriction decreases each of the positive and negative
eigenvalue counts by at most one. Since
\[
 \nu_+\!\left(q_\ell\big|_{x_{-\lambda}=0}\right)
 +\nu_-\!\left(q_\ell\big|_{x_{-\lambda}=0}\right)
 =p+r-2=\nu_+(q_\ell)+\nu_-(q_\ell)-2,
\]
both counts decrease by one.
Thus, for $\lambda\le\ell\le n_0$,
\[
 \nu_-\!\left(q_\ell\big|_{x_{-\lambda}=0}\right)
 =p+r-\ell-\mathbf1_{\{n_0+\ell\text{ odd}\}},\qquad
 \nu_0\!\left(q_\ell\big|_{x_{-\lambda}=0}\right)
 =2\mathbf1_{\{n_0+\ell\text{ odd}\}}.
\]
For $\ell>n_0$, restriction removes one positive eigenvalue.

In the original chain coordinates, let $x$ satisfy only the constraint
from $VZ^t=0$, imposed when $\ell\ge\lambda$.
Let $u$ be its Euclidean orthogonal projection onto the coefficient space
of $M_\ell$. Since $a_{\lambda,n_0}^2+a_{\lambda+2,n_0}^2=4/K$,
the projection is given by $u_0^{(1)}=0$ and
\[
 \begin{pmatrix}u_{-1}^{(1)}\\u_1^{(1)}\end{pmatrix}
 =\begin{pmatrix}x_{-1}^{(1)}\\x_1^{(1)}\end{pmatrix}
 -\frac K4\bigl(a_{\lambda,n_0}x_{-1}^{(1)}
                 -a_{\lambda+2,n_0}x_1^{(1)}\bigr)
 \begin{pmatrix}a_{\lambda,n_0}\\-a_{\lambda+2,n_0}\end{pmatrix},
\]
with $u_q^{(i)}=x_q^{(i)}$ for all remaining coordinates.
Thus only the coordinates $0,\pm1$ in the first chain change,
so the constraint from $VZ^t=0$ is preserved. Completing squares gives
\begin{equation}\label{direct:nonzero-squares}
\begin{aligned}
 \sum_{i=1}^{c}q_\ell(x^{(i)})
 ={}&\bigl(\sqrt L\,x_0^{(1)}
       -a_{\lambda,n_0}x_{-1}^{(1)}
       -a_{\lambda+2,n_0}x_1^{(1)}\bigr)^2\\
 &+\frac K4\bigl[\sqrt{L-2}
       \bigl(a_{\lambda,n_0}x_{-1}^{(1)}
             -a_{\lambda+2,n_0}x_1^{(1)}\bigr)
       -a_{\lambda,n_0}a_{\lambda-2,n_0}x_{-2}^{(1)}
       +a_{\lambda+2,n_0}a_{\lambda+4,n_0}x_2^{(1)}\bigr]^2
       +M_\ell(u).
\end{aligned}
\end{equation}
For $\ell\ge2$, the two displayed linear forms and the orthonormal
coordinates of $u$ in this space give an invertible
triangular change of coordinates, since $\sqrt L,\sqrt{L-2}>0$.
At $\ell=1$, the second square vanishes, while
\[
 a_{\lambda,n_0}x_{-1}^{(1)}-a_{\lambda+2,n_0}x_1^{(1)}
\]
is retained as a coordinate on which the quadratic form does not depend.
Together with the linear form in the first square of
\eqref{direct:nonzero-squares} and the coordinates of $u$,
it gives an invertible change of coordinates.
On the coefficient space subject to $VZ^t=0$, the left-hand side
of \eqref{direct:nonzero-squares} is therefore congruent to
\[
 \operatorname{diag}(1,1,M_\ell)\quad(\ell\ge2),\qquad
 \operatorname{diag}(1,0,M_1)\quad(\ell=1).
\]

For $1\le\ell\le n_0$, the inertia counts for the restricted forms and
\eqref{direct:nonzero-squares} give
\[
 \nu_-(M_\ell)
 =c(p+r-\ell)-\mathbf1_{\{\ell\ge\lambda,\ n_0+\ell\text{ odd}\}},\qquad
 \nu_0(M_\ell)
 =c-\mathbf1_{\{\ell\ge\lambda\}}
   +2\mathbf1_{\{\ell\ge\lambda,\ n_0+\ell\text{ odd}\}}
   -\mathbf1_{\{\ell=1\}}.
\]
For $\ell>n_0$, the restricted scalar form is positive definite,
so the same congruence gives $M_\ell>0$.
Since $\nu_+(J_\ell)=2\nu_-(M_\ell)$ and
$\nu_0(J_\ell)=2\nu_0(M_\ell)$, these counts give
Proposition~\ref{index:inertia}\textup{(i)}.

\subsection{Proof of Proposition~\ref{index:inertia}\textup{(ii)}}\label{app:zero-congruence}

To determine the inertia of $A_\ell$, note that for $\ell\ge2$ we have
$v^tTv=4(L-2)/K$, and multiplication by $T$ gives
\[
 c^tTv=-\sqrt{L-2}
 \left(\sqrt{\alpha_1(\alpha_1-2)}\,c_{1,2}
       -\sqrt{\alpha_2(\alpha_2-2)}\,c_{2,2}\right).
\]
Writing $x=c+tv$ with $c\perp v$ and completing the square gives
\[
 x^tTx=A_\ell(c)+\frac{4(L-2)}K
       \left(t+\frac{Kc^tTv}{4(L-2)}\right)^2.
\]
Therefore $T$ is congruent to $\operatorname{diag}(L-2,A_\ell)$.
For $n_0=m_i$ and $\lambda=0$, the scalar difference operator in
Appendix~\ref{app:nonzero-congruence} preserves parity and $\omega_{-q}=\omega_q$.
Its eigenpolynomial of degree $j$ has parity $j$. Since $T_i$
represents the restriction to $u_{-q}=-u_q$,
\eqref{direct:scalar-spectrum} gives its eigenvalues as
\[
 (\ell-j)(m_i+\ell-j+1),\qquad j=1,3,\ldots,2d_i-1.
\]
Their signs are those of $\ell-j$.
For $\ell=1$, $T=0$ and $A_1=0$ on $v^\perp$.
For all $\ell\ge1$, counting the odd indices therefore gives
\[
 \nu_-(A_\ell)
 =\sum_{i=1}^2\max\left\{0,d_i-\left\lfloor\frac{\ell+1}{2}\right\rfloor\right\},\qquad
 \nu_0(A_\ell)
 =\mathbf1_{\{\ell\text{ odd}\}}\sum_{i=1}^2\mathbf1_{\{\ell<m_i\}}
   -\mathbf1_{\{\ell=1\}}.
\]

We set $\delta=(m_1-m_2)/2$. The eigenvalues of $D$ are
\[
 e_{i,j}=(\ell-j)(m_i+\ell-j+1)+(-1)^{i-1}(\alpha_1-\alpha_2)
        =(\ell+(-1)^{i-1}\delta-j)(\ell+n/2-j+1),
\]
where $i=1,2$ and $j=1,3,\ldots,2d_i-1$.
For $1\le\ell\le n/2$, we pair an odd index $j$ in the first chain
with $2\ell-j$ in the second whenever both occur. Then
\[
 e_{1,j}=(\ell+\delta-j)(\ell+n/2-j+1),\qquad
 e_{2,2\ell-j}=-(\ell+\delta-j)(n/2-\ell+j+1).
\]
Both second factors are positive, and the number of pairs is
\[
 \#\{j\text{ odd}:2\ell-2d_2+1\le j\le2d_1-1\}=d_1+d_2-\ell.
\]
All unpaired eigenvalues are positive. Each nonzero pair contributes
one negative eigenvalue, while a zero pair contributes two zero
eigenvalues. At most one pair is zero. Thus, for $1\le\ell\le n/2$,
\[
 \nu_-(D)+\frac12\nu_0(D)=d_1+d_2-\ell,
 \qquad \nu_0(D)\in\{0,2\}.
\]
For $\ell>n/2$, all $e_{i,j}$ are positive, so $D>0$.

It remains to determine the inertia of $B_\ell$.
Use the same formula to define $e_{i,j}$ for even $j$. Set
\[
 P_i=\prod_{j=0,2,\ldots,2d_i}e_{i,j},\qquad
 Q_i=\prod_{j=1,3,\ldots,2d_i-1}e_{i,j},
\]
where empty products equal $1$.
For $n_0=m_i$ and $\lambda=0$, restrict the matrix of $q_\ell$ in
Appendix~\ref{app:nonzero-congruence} to $u_{-q}=u_q$.
Express this restriction in the orthonormal coordinates
$(u_0,\sqrt2u_1,\ldots,\sqrt2u_{d_i})$ and add
$(-1)^{i-1}(\alpha_1-\alpha_2)I$.
The resulting matrix has determinant $P_i$ and first diagonal entry
$L+(-1)^{i-1}(\alpha_1-\alpha_2)$; its first off-diagonal entry, when
present, is $-\sqrt{2\alpha_iL}$.
Deleting the first row and column leaves the $i$th block of $D$,
whose determinant is $Q_i$.
Expanding the two restricted determinants along their first rows and columns
gives the first equality below; row and column operations using the
bottom-right entry $2L$ give the second:
\begin{equation}\label{directIII:productidentity}
 P_1Q_2+P_2Q_1
 =\det\begin{pmatrix}
 D&\sqrt{2L}\,v\\
 \sqrt{2L}\,v^t&2L
 \end{pmatrix}
 =2L\det B_\ell.
\end{equation}

Suppose $1\le\ell\le n/2$. If $\ell\ge\delta$, the indices
$\ell+\delta$ and $\ell-\delta$ lie in $[0,2d_1]$ and $[0,2d_2]$,
respectively, and
\[
 e_{1,\ell+\delta}=e_{2,\ell-\delta}=0.
\]
These indices have the same parity, so either $P_1=P_2=0$ or $Q_1=Q_2=0$.
Hence \eqref{directIII:productidentity} gives $\det B_\ell=0$.
If $\ell<\delta$, then $d_1=\ell$ and $Q_1Q_2\ne0$.
For $d_2=m_2/2<\ell$, the equality $n/2+1-\delta=2d_2+1$ shows that
the numerator and denominator of $P_1/Q_1$ share the factor
\[
 \prod_{j=0}^{2(\ell-d_2)-1}(\ell+\delta-j).
\]
After cancelling it, shift the index by $2(\ell-d_2)$ in the remaining
first factors. These cancellations, or direct substitution when $d_2=\ell$, give
\[
 \frac{P_1}{Q_1}=
 \frac{\prod_{j=0,2,\ldots,2d_2}
       (2d_2-\ell+\delta-j)(\ell+n/2+1-j)}
      {\prod_{j=1,3,\ldots,2d_2-1}
       (2d_2-\ell+\delta-j)(\ell+n/2+1-j)}
 =-\frac{P_2}{Q_2}.
\]
In the last equality, $j\mapsto2d_2-j$ reverses the first factors.
The numerator has one more negative sign than the denominator.
Thus \eqref{directIII:productidentity} again gives $\det B_\ell=0$.

Reverse the order of the first chain and keep the second in its original
order. Then $B_\ell$ is tridiagonal with all adjacent off-diagonal entries
nonzero. Successive rows determine every remaining coordinate of a kernel
vector from its first coordinate, so $\dim\ker B_\ell\le1$.
Together with $\det B_\ell=0$, this yields $\nu_0(B_\ell)=1$.
Since $B_\ell=D-vv^t$,
\[
 \nu_+(D)-1\le\nu_+(B_\ell)\le\nu_+(D),\qquad
 \nu_-(D)\le\nu_-(B_\ell)\le\nu_-(D)+1.
\]
Moreover,
\[
 \nu_+(D)+\nu_-(D)+\nu_0(D)=d_1+d_2
 =\nu_+(B_\ell)+\nu_-(B_\ell)+1.
\]
Together with $\nu_0(D)\in\{0,2\}$, these relations give
\[
 \nu_-(B_\ell)=\nu_-(D)+\frac12\nu_0(D)
 =d_1+d_2-\ell=\min\{\ell,m_2/2,n/2-\ell\},\qquad
 \nu_0(B_\ell)=1.
\]
For $\ell>n/2$, all factors in $P_i,Q_i$ are positive, so
\eqref{directIII:productidentity} gives $\det B_\ell>0$.
Since $D>0$ and subtracting $vv^t$ creates at most one negative
eigenvalue, it follows that $B_\ell>0$.
Together with the counts for $A_\ell$ and \eqref{directIII:J},
this proves Proposition~\ref{index:inertia}\textup{(ii)}.

\section*{Acknowledgements}
The authors thank Professors Xiaowei Xu and Hongbin Cui for their help
and encouragement.

\section*{Funding}
This work is supported by the National Natural Science Foundation of China
(Grant No.\ 12571057), the project of Stable Support for Youth Team in Basic
Research Field of the Chinese Academy of Sciences (Grant No.\ YSBR-001),
and the Fundamental Research Funds for the Central Universities
(Grant No.\ WK0010000081).

\section*{Declarations}
\paragraph{Author contributions}
Qing Chen: Writing -- review \& editing (equal).
Shuping Huang: Writing -- review \& editing (equal).
Shoulan Yi: Writing -- review \& editing (equal).
\paragraph{Conflict of interest}
The authors have no conflicts of interest to disclose.
\paragraph{Data availability}
Data sharing is not applicable to this article, as no datasets were generated
or analysed.
\paragraph{Use of generative AI}
The authors used ChatGPT (OpenAI) to assist with language editing.
The authors reviewed and edited the content and take full responsibility
for the final manuscript.

\end{document}